\documentclass[12pt]{paper}

\usepackage{amssymb,amsfonts,amsthm,amsmath,comment}

\usepackage{latexsym,multicol,graphicx}

\usepackage{verbatim}
\usepackage{xcolor}
\usepackage{verbatim}
\DeclareGraphicsExtensions{.pdf,.png,.jpg}
\usepackage{pict2e}

\usepackage[a4paper, hmargin=2.5cm, vmargin={3cm, 3cm}]{geometry}

\usepackage{fancyhdr}
\usepackage[title,titletoc,toc]{appendix}

\newtheorem{theorem}{Theorem}
\newtheorem{lemma}{Lemma}

\newtheorem{claim}{Claim}

\renewcommand{\phi}{\varphi}

\newcommand{\e}{\varepsilon}

\newcommand{\la}{\lambda}

\newcommand{\cE}{\mathcal E}

\newcommand{\cM}{\mathcal M}
\newcommand{\cQ}{\mathcal Q}
\newcommand{\cS}{\mathcal S}

\newcommand{\bE}{\mathbb E}

\newcommand{\osc}{\operatorname{osc}}

\def\done{{1\hskip-2.5pt{\rm l}}}

\renewcommand{\le}{\leqslant}
\renewcommand{\ge}{\geqslant}

\newcommand{\bR}{\mathbb R}
\newcommand{\bC}{\mathbb C}

\newcommand{\bD}{\mathbb D}

\newcommand{\bN}{\mathbb N}

\begin{document}

\title{Translation-invariant probability measures on entire functions}

\author{Lev Buhovsky\thanks{Supported in part by ISF Grant 1380/13,
and by the Alon Fellowship},
Adi Gl\"ucksam\thanks{Supported in part by ERC Advanced Grant~692616 and ISF Grant~382/15},
Alexander Logunov\thanks{Supported in part by ERC Advanced Grant~692616 and ISF Grants~1380/13, 382/15},
Mikhail Sodin\thanks{Supported in part by ERC Advanced Grant~692616 and ISF Grant~382/15}}

\maketitle


\begin{abstract}
We study non-trivial translation-invariant probability measures on
the space of entire functions of one complex variable. The existence (and even an
abundance) of such measures was proven by Benjamin Weiss.
Answering Weiss' question, we find a relatively sharp lower bound for the
growth of entire functions in the support of such measures. The proof of
this result consists of two independent parts: the proof of the lower bound
and the construction, which yields its sharpness. Each of these parts combines
various tools (both classical and new) from the theory of entire and subharmonic
functions and from the ergodic theory.

We also prove several companion results, which concern
the decay of the tails of non-trivial trans\-la\-tion-invariant probability
measures on the space of entire functions and
the growth of locally uniformly recurrent entire and meromorphic functions.

\end{abstract}

\section{Introduction and main results}

Our starting point is Benjamin Weiss' work~\cite{Weiss} where he showed that
there exist non-trivial translation-invariant probability measures on the space
of entire functions of one complex variable (a formal definition of such measures
will be given several lines below). Actually, Weiss showed that there is an abundance
of such measures. Another approach to the construction of such measures was suggested
by Tsirelson~\cite{Tsirel}. Tsirelson worked in a somewhat different and simpler
context. The works of Weiss and Tsirelson raise a number of intriguing questions which lie
at the crossroads of complex analysis and ergodic theory. Here, we address some of
these questions.

\subsection{}
Let $\cE$ denote the space of entire functions with the topology generated by the
semi-norms
\[
\| F \|_K = \max_K | F |\,
\]
where $K$ runs over all compact subsets of $\bC$, and let $B$ be the Borel sigma-algebra
generated by this topology. Then $\bC$ acts on $(\cE, B)$ by translations:
\[
(\tau_w F)(z) = F(z+w), \quad w\in\bC\,.
\]
The probability measure $\la$ on $(\cE, B)$ is called {\em translation-invariant}
if it is invariant with respect to this action. The translation-invariant measure $\la$ is
called {\em non-trivial} if the set of all constant functions in $\cE$ has measure zero
(the constant functions are fixed points of the action $\tau$). Due to~\cite{Weiss}, we know
that non-trivial translation invariant probability measures on $(\cE, B)$ exist.
In what follows, we retain the notation $\la$ for such measures.

After some reflection it becomes plausible that entire functions from the Borel support of $\la$ must grow sufficiently fast and that $\la$ must have heavy tails. The goal of this work is to
justify these statements.

\subsection{}
For an entire function $F$ we put
\[
M_F(R) = \max_{R\,\overline\bD} |F|\,,
\]
where $R\,\overline\bD=\{z\colon |z|\le R\}$.

\begin{theorem}\label{thm1} \mbox{}

\smallskip\noindent {\rm (A)} Let $\la$ be a non-trivial translation-invariant
probability measure on the space of entire functions. Then, for
$\la$-a.e. function $F$ and for every $\e>0$,
\[
\lim_{R\to\infty} \frac{\log\log M_F(R)}{\log^{2-\e} R} = +\infty\,.
\]

\smallskip\noindent {\rm (B)}
There exists a non-trivial translation-invariant probability measure on
the space of entire functions such that, for $\la$-a.e. function $F$ and
for every $\e>0$,
\[
\lim_{R\to\infty} \frac{\log\log M_F(R)}{\log^{2+\e} R} = 0\,.
\]
\end{theorem}

\subsubsection{}\label{subsub-1.2.1}

The proof of the first part of Theorem~\ref{thm1} relies on a growth estimate of subharmonic
functions, which might be of independent interest. To bring this estimate
we introduce some notation.
\begin{itemize}
\item
Till the end of Section~\ref{subsub-1.2.1}, we
assume that {\em all squares denoted by $Q$ and $S$ have all four vertices
with integer-valued coordinates}.
\end{itemize}

Let $u$ be a non-negative subharmonic function on a neighbourhood of the
square $Q\subset \bC$ with side-length $L(Q)$. Let
$ M_u(Q)=\max_{\bar Q} u $ and $ Z_u = \{ u=0 \} $. We denote by $A$ the
area measure and by $|X|$ the cardinality of a finite set $X$.

Given $\gamma\in (0, 1)$, we say that a unit square $S$ (i.e., the square with
$L(S)=1$) is $\gamma$-{\em good} if (i) $A(S\cap Z_u)\ge\gamma$ and (ii) $M_u(S)~\ge~1$.
For any square $Q$, we put
\[
\beta (Q) = \beta_{u, \gamma} (Q) = \frac{\bigl| \{S\subset Q\colon S\ {\rm is\ } \gamma{\rm-good\ unit\ square} \} \bigr|}{A(Q)}\,.
\]

\begin{lemma}\label{lemma-LevSasha}
Given $\gamma, \beta\in (0, 1)$ there exists $c=c(\gamma, \beta)>0$
such that for any square $Q$ with $L(Q)=L \ge 10$ and any non-negative
subharmonic function $u$ on a neighbourhood of $Q$ with $\beta=\beta(Q)$,
\[
M_u(Q) \ge e^{c\left( \frac{\log L}{\log\log L} \right)^2}\,.
\]
\end{lemma}

It is instructive to juxtapose this estimate with a less restrictive one
(which also will be used below), where we require that
{\em almost every} unit square $S\subset Q$ contains a non-negligible piece of the set $Z_u$
and get a much faster growth of $u$.
\begin{lemma}\label{lemma-Adi}
Let $u$ be a non-negative subharmonic function on a neighbourhood of the square $Q$
with $L=L(Q)$ and let $\alpha>0$ be a positive parameter.
Suppose that for some $\gamma\in (0, 1)$ and for all, except of at most $\alpha L$
unit squares $S\subset Q$, we have $A(S\cap Z_u) \ge \gamma$. Then,
\[
M_u(Q) \ge e^{cL} M_u\bigl([-\tfrac12, \tfrac12]^2\bigr)
\]
with some $c=c(\gamma, \alpha)>0$, provided that the size $L$ of the square $Q$ is sufficiently large.
\end{lemma}
Note that our reduction of the first part
of Theorem~\ref{thm1} to Lemma~\ref{lemma-LevSasha} is based on the
pointwise ergodic theorem, and that the proof of Lemma~\ref{lemma-LevSasha}
makes use of Lemma~\ref{lemma-Adi}.

\subsubsection{}

A natural idea to construct a non-trivial translation-invariant probability
measure on $\cE$ (and, in particular, for the proof of the second part of Theorem~\ref{thm1}) is to use the classical Krylov-Bogolyubov construction. We take a function $F\in\cE$,
denote by $\delta_F$ the point mass on $F$ (viewed as a probability measure on $\cE$)
and average it along the orbit of $\tau$ defining
\[
\la_R = \frac1{\pi R^2}\, \int_{R\,\bD} \delta_{\tau_w F}\, {\rm d}A(w)\,, \quad R\ge 1\,.
\]
In other words, for any Borel set $B\subset \cE$,
\[
\la_R (B) = \frac1{\pi R^2}\, \int_{R\,\bD} \done_B(\tau_w F)\, {\rm d}A(w)\,.
\]

Then, we let $R\to\infty$, and consider the limiting measure.
The problem with this idea is that the space $\cE$ is not compact;
therefore, we need to ensure tightness of the family $(\la_R)_{R\ge 1}$.
In addition, we must ensure that (at least a part of) the limiting measure
is not supported by the constant functions. Thus, the entire function $F$, which
we start with, should be carefully chosen.

First, we construct a particular subharmonic function $u$ which can be thought
as a certain approximation to $\log |F|$. We define a special
unbounded closed set $E\subset \bC$ which can be thought as
a two-dimensional fat Cantor-type set viewed from the inside-out and a subharmonic
function $u$ of a nearly minimal growth outside $E$ (Lemma~\ref{lemma-SH}).
Then, using H\"ormander's classical estimates of solutions to $\bar\partial$-equations,
we build an entire function $G$ of a nearly minimal growth outside $E$
with needed properties (Lemma~\ref{lemma-EF}). The functions $u$ and $G$ enjoy
an interesting dynamical behaviour, and their construction is likely of
independent interest.

\subsection{}

We say that an entire function $F$ is {\em locally uniformly recurrent} if for
every $\e>0$ and every compact set $K\subset\bC$ the set
$ \bigl\{w\in\bC\colon \max_K |\tau_w F - F | < \e \bigr\} $
is {\em relatively dense} in $\bC$ (that is, any disk of sufficiently large radius contains
at least one point of this set). This is a locally uniform counterpart of Bohr's classical
definition of almost-periodicity. In~\cite{Weiss}, Weiss gave a simple construction of functions of this class based on the Runge approximation theorem.

Locally uniformly recurrent entire functions can serve as a starting point for the
Krylov-Bogolyubov-type construction described above. However, as the following
theorem shows their growth is rather far from the minimal one.

\begin{theorem}\label{thm-unform-recurrent}\mbox{}

\smallskip\par\noindent{\rm (A)} For any non-constant locally uniformly recurrent entire
function $F$,
\[
\liminf_{R\to\infty} \frac{\log\log M_F(R)}{R} > 0\,.
\]

\smallskip\par\noindent{\rm (B)} There exists a non-constant locally uniformly recurrent
functions $F$ such that
\[
\limsup_{R\to\infty} \frac{\log\log M_F(R)}{R} < \infty\,.
\]
\end{theorem}

Note that the difference in the growth of entire functions in Theorems~\ref{thm1} and~\ref{thm-unform-recurrent} and that of the corresponding
subharmonic functions in Lemmas~\ref{lemma-LevSasha} and~\ref{lemma-Adi} are
closely related.

\subsection{}

As we have already mentioned, translation-invariant probability measures on
the space of entire functions must have heavy tails.
\begin{theorem}\label{thm-heavytails}\mbox{}

\smallskip\par\noindent{\rm (A)} Let $\la$ be a non-trivial translation-invariant probability measure on the space of entire functions. Then, for every $\e>0$,
\[
\bE \bigl[ (\log\log |F(0)|)^{1+\e} \bigr] = +\infty\,.
\]

\smallskip\par\noindent{\rm (B)} There exists a non-trivial translation-invariant probability measure $\la$ on the space of entire functions such that, for every $t\ge 1$,
\[
\la \bigl\{F\colon \log\log |F(0)| > t \bigr\} \lesssim \frac1{t}\,.
\]
\end{theorem}

Here and elsewhere, the notation $X \lesssim Y$  means that there exists a positive numerical constant
$C$ such that $X \le CY$.

\subsection{}
It is natural to look at the counterparts of Theorems~\ref{thm1} and~\ref{thm-unform-recurrent} for meromorphic functions. We treat
meromorphic functions as maps of the complex plane into the Riemann sphere endowed
with the spherical metric $\rho$, and denote by $\cM$ the space of meromorphic
functions endowed with the topology of the locally uniform convergence in the spherical
metric (as usual, we treat $\infty$ as a constant
meromorphic function). By $B$ we denote the Borel sigma-algebra generated by this topology.
Since $\cE\subset\cM$, it is worthwhile to note that these definitions are consistent
with the ones we have used above.

To measure the growth of a
meromorphic function $F$ we will use Nevanlinna's characteristics $T_F(R)$.
It will be convenient to use it in
the Ahlfors-Shimizu geometric form:
\[
T_F(R) = \int_0^R \Bigl( \frac1{\pi}\, \int_{r\,\bD} F^{\#}(z)^2 \, {\rm d}A(w)\Bigr) \, \frac{{\rm d}r}r\,,
\]
where
\[
F^{\#}(z) = \frac{|F'(w)|}{1+|F(w)|^2}
\]
is the spherical derivative of $F$. Then the inner integral in the definition of
characteristics $T_F$ is the spherical area of the image of the disk $F(r\,\bD)$ considered with multiplicities
of covering. The basic properties of the Nevanlinna's characteristics can be found,
for instance, in~\cite[Chapter~1]{GO}. Here, we will mention that if $F$ is an
entire function then the growth of its Nevanlinna characteristics and of the logarithm
of its maximum modulus are equivalent in the following sense:
\[
T_F(R) < \log M_F(R) + O(1)\,,
\]
and, for every $R_1>R$,
\[
\log M_F(R) < \frac{R_1+R}{R_1-R}\, T_F(R_1) + O(1)\,.
\]
We also point out that it is easy to see that if $F$ is a non-constant doubly periodic meromorphic function, then
the spherical area of the image  $F(r\,\bD)$ counted with multiplicities has quadratic grows with $r$, and therefore, $T_F(R)$ has quadratic growth as well.

As above, $\bC$ acts on $\cM$ by translations.
We call the probability measure $\la$ on $(\cM, B)$ translation-invariant
if it is invariant with respect to this action. As above, we call a translation-invariant
measure $\la$ non-trivial if the set of all non-constant functions has measure zero.
Examples of non-trivial translation-invariant probability measures can be easily constructed
by averaging the translations of a doubly periodic function. In these examples, for $\la$-a.e.
function $F\in\cM$, $T_F(R)=O(R^2)$ as $R\to\infty$. The following theorem shows that one
cannot do better:
\begin{theorem}\label{thm-meromorh-measures}
Let $\la$ be a non-trivial translation-invariant
probability measure on meromorphic functions. Then, for $\la$-a.e. function $F\in\cM$,
\[
\liminf_{R\to\infty} \frac{T_F(R)}{R^2} > 0\,.
\]
\end{theorem}

We call a meromorphic function $F$ {\em locally uniformly recurrent}
if for every $\e>0$ and every
compact set $K\subset \bC$, the set
$ \bigl\{w\in\bC\colon \max_K \rho(\tau_w f, f) < \e \bigr\} $
is relatively dense in $\bC$. Here, as above, $\rho$ is the spherical distance.
It is easy to see that doubly periodic meromorphic functions are locally uniformly
recurrent. I.e., there are plenty of locally uniformly recurrent meromorphic functions
$F$ with $T_F(R)=O(R^2)$ as $R\to\infty$. As in the previous case, this estimate cannot
be improved:
\begin{theorem}\label{thm-meromorh-recurrent}
Let $F$ be a non-constant locally uniformly recurrent
meromorphic function. Then
\[
\liminf_{R\to\infty} \frac{T_F(R)}{R^2} > 0\,.
\]
\end{theorem}

\subsection*{Acknowledgments}

The authors are grateful to Alexander Borichev, Gady Kozma, Fedor Nazarov and Benjy Weiss for
several very helpful discussions. We thank Steven Britt for his assistance with
language editing.

\section{Proof of Lemmas~\ref{lemma-LevSasha} and~\ref{lemma-Adi}}

In this section the squares denoted by $Q$, $Q_j$, $\cQ$, and $S$ have vertices with integer-valued coordinates,
$Q$ is a square with large side-length $L=L(Q)$, and $u$ is a subharmonic function
on a neighbourhood of $Q$. By $Z_u=\{u=0\}$ we denote the zero set of $u$.

\subsection{Proof of Lemma~\ref{lemma-Adi}}
Assuming that for all but $\alpha L$ unit squares $S\subset Q$ we have
$A(S\cap Z_u)\ge \gamma$, we need to show that
\[
\max_Q u \ge e^{c(\alpha, \gamma)L} \max_{[-1/2, 1/2]^2} u\,.
\]

\medskip
First, we observe that if the disk $D_z\subset Q$ centered at $z$
contains a portion of the zero set $Z_u$ of area at least
$\gamma$  (with $\gamma< A(D_z)$), then
\[
u(z) \le \frac1{A(D_z)}\, \int_{D_z} u\, {\rm d}A
= \frac1{A(D_z)}\, \int_{D_z\setminus Z_u} u\, {\rm d}A
\le \bigl( 1-\frac{\gamma}{A(D_z)} \bigr)\, \max_{\bar D_z} u\,,
\]
whence,
\[
\max_{\bar D_z} u \ge \Bigl( 1-\frac{\gamma}{A(D_z)} \Bigr)^{-1} u(z)\,.
\]

Let $N$ be the integer part of $\tfrac12 L$. Put
$M_u(r)\stackrel{\rm def}=\max \bigl\{u(z)\colon |z|\le r \bigr\}$,
take the points $z_0=0$, $z_1$, \ldots , $z_N$, with $|z_j|=j$,
so that $u(z_j)=M_u(j)$, $j=1, \ldots , N$,
and consider the disks $D_j=D(z_j, p)$ with sufficiently large integer $p\ge 2$.
We call the index $j\le N-p$ {\em normal} if the disk $D_j$ contains at least one non-exceptional unit square $S$.
For normal indices $j$, we have
\begin{equation}\label{eq-max}
M_u(j+p) \ge \max_{\bar D_j} u \ge \Bigl( 1 - \frac{\gamma}{\pi p^2} \Bigr)^{-1} u(z_j)
= \Bigl( 1 - \frac{\gamma}{\pi p^2} \Bigr)^{-1} M_u(j)\,.
\end{equation}

If the disks $D_{j_1}$, \ldots , $D_{j_\ell}$ are not normal, then
the number of different exceptional squares contained in their union $ D_{j_1} \bigcup \ldots \bigcup D_{j_\ell}$ is $\gtrsim \ell p$. Since the total number of exceptional unit squares does not exceed $\alpha L$, we conclude that the number of not normal disks is
$\lesssim \alpha p^{-1} L < \tfrac15 L$ provided that $p$ was chosen much larger than $\alpha$. We conclude that there are at least $\tfrac14 L$ indices $1\le j \le N-p$,
for which estimate~\eqref{eq-max} holds. Hence, the lemma follows. \hfill $\Box$

\subsection{Proof of Lemma~\ref{lemma-LevSasha}}
Recall that we say that a unit square $S\subset Q$ is $\gamma$-good if
$A(S\cap Z_u)\ge \gamma$ and $\max_S u \ge 1$, and that for any square $Q$,
we put
\[
\beta (Q) = \frac{\bigl| \{S\subset Q\colon S\ {\rm is\ } \gamma{\rm-good\ unit\ square} \} \bigr|}{A(Q)}\,.
\]
Our aim is to show that
\[
\log \max_Q u \ge c(\beta, \gamma) \Bigl( \frac{\log L}{\log\log L} \Bigr)^2,
\quad \beta=\beta(Q)\,.
\]

\medskip
With no loss of generality we assume that $L=k^k$ for an integer $k$
(then, $k\simeq \tfrac{\log L}{\log\log L}$). We
construct a sequence of squares $Q_0=Q$, \ldots , $Q_k$, with $L(Q_j)=k^{k-j}$.
First, we split the square $Q_{j-1}$ into $k^2$ squares $\cQ$ with
$L(\cQ)=k^{k-j}$. For these squares $\cQ$ we write $\cQ\prec Q_{j-1}$, and note that
\begin{equation}\label{eq-beta}
\beta (Q_{j-1}) = \frac1{k^2}\, \sum_{\cQ\prec Q_{j-1}} \beta (\cQ)\,.
\end{equation}
Then, according to certain rules described below, we choose one of the squares $\cQ$,
and call it $Q_j$.

Suppose that the squares $Q_0$, \ldots , $Q_{j-1}$ have already been chosen.
We will fix the parameters $B>1$ and $0<\theta<1$ to be chosen later, and consider three cases.

\medskip\noindent\underline{Case 1}: {\em there exist at least $B k$ squares $\cQ\prec Q_{j-1}$ such that
$\beta (\cQ) < \tfrac12 \beta (Q_{j-1})$}.

\medskip\noindent
We claim that in this case there exists at least one square $\cQ\prec Q_{j-1}$ with
\begin{equation}\label{eq-beta'}
\beta (\cQ) \ge \Bigl( 1+\frac{B}{2k} \Bigr) \beta(Q_{j-1})\,.
\end{equation}
Indeed, otherwise, \eqref{eq-beta} gives us
\[
1 \le \frac1{k^2}
\Bigl( Bk \cdot \frac12 + (k^2-Bk) \cdot \Bigl( 1 + \frac{B}{2k} \Bigr) \Bigr)
= \frac1{k^2} \bigl( k^2 - \tfrac12 B^2 \bigr) < 1\,,
\]
arriving at a contradiction.

Then, we let $Q_j$ be one of the squares $\cQ\prec Q_{j-1}$ such that~\eqref{eq-beta'}
holds.

\medskip\noindent\underline{Case 2}: {\em for all squares
$\cQ\prec Q_{j-1}$ contained in the square $(1-\theta) Q_{j-1}$},
\[
\beta (\cQ) < \bigl( 1- \frac{1}k \bigr) \beta (Q_{j-1})\,.
\]
Here, $(1-\theta)Q_{j-1}$  denotes the square with the same center as $Q_{j-1}$
and $L((1-\theta)Q_{j-1})=(1-\theta)L(Q_{j-1})$.

\medskip\noindent
We claim that if $\theta$ is chosen sufficiently small, then~\eqref{eq-beta'} holds
for at least one of the remaining squares. Otherwise,
\begin{multline*}
1 < \frac1{k^2} \Bigl( (1-\theta)^2k^2 \cdot \Bigl(1-\frac1{k} \Bigr)
+ \bigl( 1 - (1-\theta)^2 \bigr)k^2 \cdot \Bigl(1+\frac{B}{2k} \Bigr) \Bigr) \\
= \frac1{k^2}
\Bigl( k^2 - k \bigl( (1-\theta)^2 - \tfrac12 B(1-(1-\theta)^2 \bigr) \Bigr)
< 1 - \frac1{k} \bigl( (1-\theta)^2 - B\theta \bigr) < 1\,,
\end{multline*}
provided that $B\theta < \tfrac12$.

As in the first case, we let $Q_j$ be one of the squares $\cQ\prec Q_{j-1}$
such that~\eqref{eq-beta'} holds.

\medskip
We now consider the remaining case, which is
complementary to the cases 1 and 2:

\medskip\noindent\underline{Case 3}: {\em
there exists at least one square $\cQ\prec Q_{j-1}$  contained in $(1-\theta)Q_{j-1}$
such that
$\beta (\cQ) \ge \bigl(1-\tfrac1{k}\bigr)\beta(Q_{j-1})$ (with $\theta=\theta(B)$ chosen above)}.
At the same time, the number of squares $Q\prec Q_{j-1}$  with $\beta (Q) \ge \frac12 \beta(Q_{j-1})$ is
not less than $k^2-Bk$.

\medskip\noindent
Then we
call one one these squares $Q_j$. We also know that {\em for at most $Bk$ squares  $Q\prec Q_{j-1}$,
we have  $\beta(Q)\le \tfrac12 \beta(Q_{j-1})$}.

\medskip Now, we are ready to prove Lemma~\ref{lemma-LevSasha}. First, we note that
on each step the value $\beta (Q_j)$ either increases (cases 1 and 2), or decreases
by a factor of at most $1-\tfrac1{k}$.
Since the total number of steps is $k\ge 2$, we conclude that {\em for each $j$, $\beta (Q_j)\ge  \tfrac13 \beta (Q_0) =
\tfrac13 \beta$}.

\medskip
Next, we observe that if
on the $j$th step one of the cases 1 or 2 occurs, then  by~\eqref{eq-beta'},
$\beta(Q_j)$ will increase by a factor of at least $1+(2k)^{-1} B$.
Since on other steps $\beta (Q_j)$
decreases not more than a factor of $1-\tfrac1{k}$, choosing $B=B(\beta) $
sufficiently large, ensures us that out of the $k$ steps at least $k/2$ steps
result in case $3$.
\begin{figure}[h!]
\centering
\includegraphics[scale=0.75]{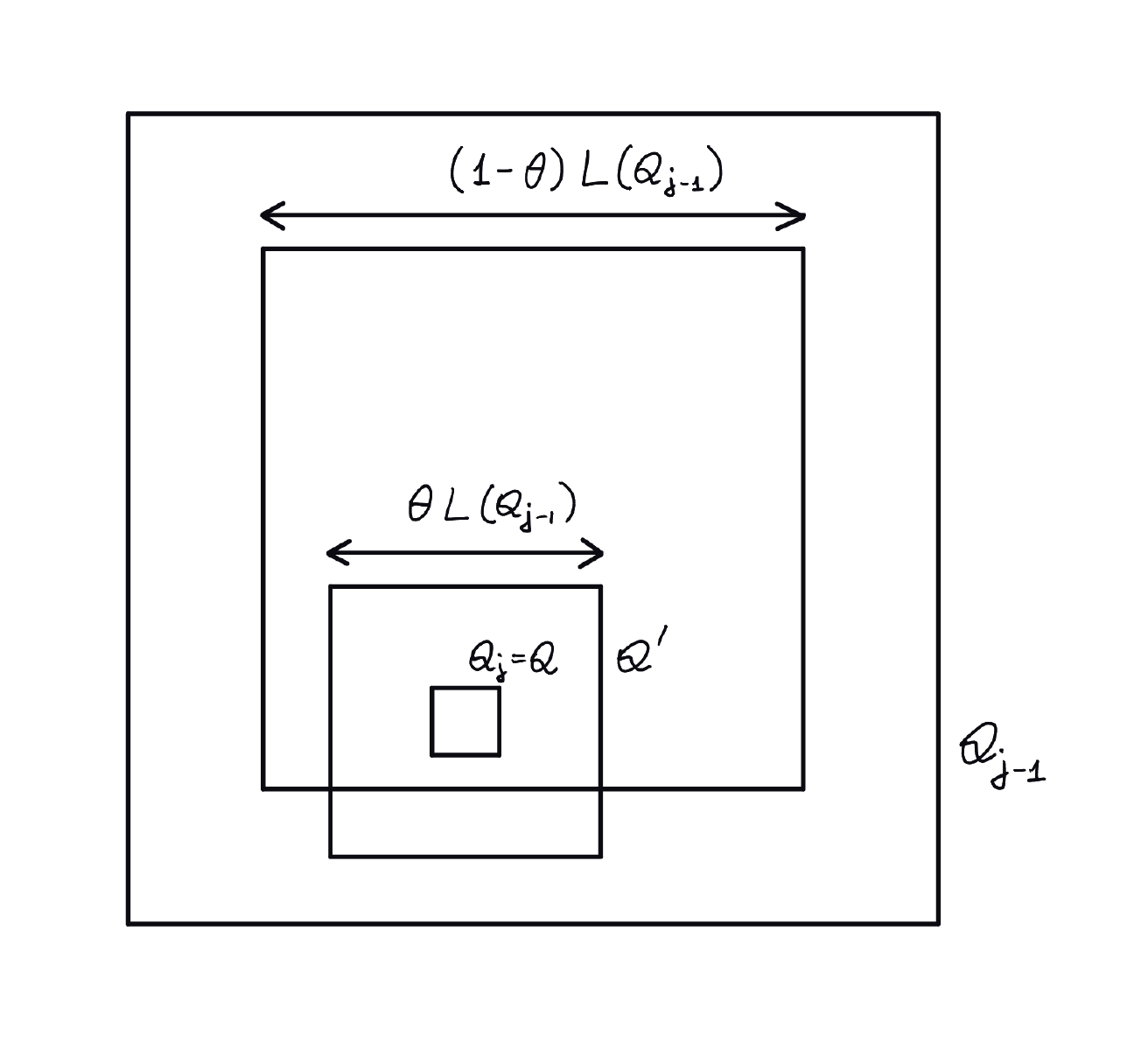}
\caption{\label{fig:case3}}
The squares $Q_{j-1}$ and $Q_j$ in the 3rd case.
\end{figure}
Assume that on the $j$th step the 3rd case happens. Then, applying Lemma~\ref{lemma-Adi} (with an
appropriate scaling) to the square
$Q'$  with $L(Q')=\theta L(Q_{j-1})$ centered at the same point as $Q_j$ (and therefore,
contained in $Q_{j-1}$),  we obtain
\[
M_u(Q_{j-1}) \ge M_u(Q') \ge e^{ck} M_u(Q_j)
\]
with some $c=c(\gamma, \beta)$.
Since this happens for at least $k/2$ indices $j$, we conclude that
\[
M_u(Q_0) \ge e^{ck^2} M_u(Q_k)\,.
\]
It remains to recall that $k\simeq \tfrac{\log L}{\log\log L}$ and that, since $Q$
was a good square, $M_u(Q_k)\ge 1$. \hfill $\Box$

\section{Proof of Theorems~1A and~2A}

In this section, $S_\rho(z)$ denotes the square of side-length $\rho$ centered at $z$,
and we let $S_\rho = S_\rho (0)$.

\subsection{An integral-geometric lemma}

We will be using a simple and known fact from the integral geometry:
\begin{lemma}\label{lemma-IG}
For any measurable set $X\subset\bC$ and any $0<\rho<R$,
\[
\Bigl|\, \frac{A(S_R\cap X)}{A(S_R)} - \frac1{A(S_R)}\,\int_{S_R}
\frac{A(S_\rho(z)\cap X)}{A(S_\rho)}\, {\rm d}A(z)\, \Bigr| \lesssim \frac{\rho}R\,.
\]
\end{lemma}

\subsection{Proof of Theorem~1A}

\subsubsection{}
Applying the ergodic decomposition theorem
(see, for instance,~\cite[Sections~6.1 and~8.6]{EW}),
we can find a Borel probabity space $(\Omega, \mathcal F, \nu)$
and a measurable map $\omega\mapsto \la_\omega$ for which

\medskip\noindent (i) for $\nu$-a.e. $\omega$,
$\la_\omega$ is a probability measure on $(\cE, B)$, which is invariant and
ergodic with respect to the action of $\bC$  on $(\cE, B)$  by translations $\tau$;

\medskip\noindent (ii) for every Borel set $\mathcal X\in B$,
\[
\la (\mathcal X) = \int_\Omega \la_\omega(\mathcal X)\, {\rm d}\nu (\omega)\,.
\]

\medskip It is not difficult to see that the set of entire functions $F$ such that, for every $\e>0$,
\[
\lim_{R\to 0}\, \frac{\log\log M_F(R)}{\log^{2-\e}R} =+\infty
\]
is a Borel set. Hence, proving Theorem~1A,
it suffices to assume that the measure $\la$ is {\em ergodic} with respect to translations $\tau$.

\subsubsection{}
Put $X(F)=\{z\in\bC\colon |F(z)|\le 1 \}$.
Given $\rho>1$, consider the Borel sets
\[
\cE_1(\rho) = \{F\in\cE\colon A(S_\rho\cap X(F)) \ge 1 \}\,,
\quad
\cE_2(\rho) = \{F\in\cE\colon \max_{\bar S_\rho} |F| \ge e \}\,.
\]
For $\rho<\rho'$, we have $\cE_i(\rho)\subset \cE_i(\rho')$, $i=1, 2$.
We denote by $\cE_i(\infty)$, $i=1, 2$, the corresponding limiting sets
as $\rho\to\infty$.
Since the complement $\cE\setminus \cE_2(\infty)$ consists of constant functions
and the measure $\la$ does not charge constants, $\la (\cE_2(\infty))=1$.

\subsubsection{}
We claim that $\la (\cE_1(\infty))=1$ as well. Otherwise, by translation-invariance of the
set $ \cE_1(\infty) $ and by ergodicity of $\la$, we have $\la (\cE_1(\infty))=0$.
Applying Lemma~\ref{lemma-IG}, we conclude that for $\la$-a.e. $F$ and for every $\rho>0$,
\[
\lim_{R\to\infty} \frac1{A(S_R)}\, \int_{S_R} \frac{A(S_{\rho}(z)\cap X(F))}{A(S_{\rho})}\,
{\rm d} A(z) = 0\,.
\]
Then, by the Wiener version of the pointwise ergodic theorem (see
for instance,~\cite{Becker}) for every $\rho>1$,
\[
\int_{\cE} A(S_\rho \cap X(F))\, {\rm d}\la (F) = 0\,,
\]
i.e., for $\la$-a.e. $F\in\cE$, $X(F)=\emptyset$, which means that $|F|\ge 1$ everywhere
in $\bC$. That is, $F$ a constant function, which is a contradiction.

\subsubsection{}
Now, we fix $\rho>1$ so that $\la \bigl( \cE_1(\rho) \cap \cE_2(\rho) \bigr) \ge \frac9{10}$,
and let
\[
X(F, \rho) = \{z\in\bC\colon A(S_\rho (z)\cap X(F))\ge 1, \, \max_{S_\rho (z)} |F| \ge e\}\,.
\]
We claim that {\em for $\la$-a.e. $F\in\cE$, the limit
\[
\lim_{R\to\infty} \frac{A(S_R\cap X(F, \rho))}{A(S_R)}
\]
exists and is $\ge\frac9{10}$}. Indeed, for any $F\in\cE$ and any $r>1$,
by Lemma~\ref{lemma-IG}, we have
\[
\frac{A(S_R\cap X(F, \rho))}{A(S_R)} =
\frac{1}{A(S_R)}\, \int_{S_R} \frac{A(S_r (z)\cap X(F, \rho))}{A(S_r)}\,
{\rm d}A(z) +O\bigl( \frac{r}{R}\bigr) \,,
\]
and by the pointwise ergodic theorem, for $\la$-a.e. $F$, the $R\to\infty$ limit
of the RHS exists and equals
\[
\int_\cE \frac{A(S_r\cap X(F, \rho))}{A(S_r)}\, {\rm d}\la (F)\,.
\]
Applying Fubini's theorem and then using the translation-invariance of the measure
$\la$, we can rewrite this expression as
\begin{align*}
& \frac1{A(S_r)}\,\int_{S_r} \Bigl[\, \int_\cE
\done_{X(F, \rho)}(z)\, {\rm d}\la (F) \,\Bigr]\, {\rm d}A(z) \\  \\
&\qquad = \frac1{A(S_r)}\,\int_{S_r} \Bigl[\, \int_\cE
\done_{X(\tau_z F, \rho)}(0)\, {\rm d}\la (F) \,\Bigr]\, {\rm d}A(z) \\ \\
&\qquad = \la\bigl\{F\in\cE\colon 0\in X(F, \rho) \bigr\} \\ \\
&\qquad = \la\bigl( \cE_1(\rho) \cap \cE_2(\rho) \bigr) \ge \frac9{10}\,,
\end{align*}
proving the claim.

\subsubsection{}
It remains to show that if $F$ is a non-constant entire function such that for some
$\rho>1$,
\[
\liminf_{R\to\infty} \frac{A(S_R\cap X(F, \rho))}{A(S_R)} \ge \frac9{10}\,,
\]
then, for every $\e>0$,
\begin{equation}\label{eq:growth}
\lim_{R\to\infty} \frac{\log\log M_F(R)}{\log^{2-\e}R} = +\infty\,.
\end{equation}

First, we note that it suffices to show that~\eqref{eq:growth} holds for the sequence
$R_n=(2\rho)^n$; then the general case follows.

Then, we take $R=(2\rho)^n$ with sufficiently large $n$, split the square $S_R$
into $ R^2/(2\rho)^2 $ squares $S$ squares $\cS$ with side-length $2\rho$, and
consider the subharmonic function $u=\log_+ |F|$. By the last claim, for
at least half of the squares $\cS$, $A(\cS\cap Z_u)\ge 1$ and $\max_{\bar\cS} u \ge 1$.
Applying Lemma~\ref{lemma-LevSasha}, we complete the proof. \hfill $\Box$

\subsubsection{Remark}

Note that with a little effort one can extract from
Lemma~\ref{lemma-LevSasha} slightly more than Theorem~1A asserts,
namely, that for $\la$-a.e. $F\in\cE$,
\[
\liminf_{R\to\infty}\, \log\log M_F(R)
\cdot \Bigl( \frac{\log\log R}{\log R} \Bigr)^2 >0\,.
\]
Likely, this estimate can be somewhat improved.

\subsection{Proof of Theorem~2A}

The proof is straightforward. Let $F$ be a non-constant locally uniformly recurrent function, and let $M=\max_{[0, 1]^2}|F|$. Applying the definition of locally uniform recurrency with
$K=[0, 1]^2$ and $\e=1$, we see that there exists $L=L(M)$ such that for every square $Q$ with the side-length $L$, $A\bigl(Q\cap \{|F|\le M+1\}\bigr) \ge 1$. Then, Lemma~\ref{lemma-Adi}
does the job. \hfill $\Box$

\section{Proof of Theorem~3A}

\subsection{A loglog-lemma that yields Theorem~3A}

We will use a version of the classical Carleman-Levinson-Sj\"oberg
loglog-theorem, cf.~\cite{Carleman, Domar0, Domar}. Likely, this lemma can be deduced from
at least one of many known versions of the loglog-theorem. Since its
proof is quite simple, for the reader's convenience, we will supply it.

\begin{lemma}\label{lemma-loglog}
Suppose $u$ is a non-constant subharmonic function in $\bC$. Then, for every
$\e>0$,
\[
\lim_{R\to\infty} \frac1{A(R\,\bD)}\, \int_{R\,\bD} \bigl( \log_+ u \bigr)^{1+\e}\, {\rm d}A
= \infty\,. \\
\]
\end{lemma}

This lemma immediately yields Theorem~3A: applying, as above,
the ergodic decomposition theorem,
we may assume that the measure $\la$ is ergodic. Then, by the pointwise ergodic theorem,
for $\la$-a.e. entire function $F$ we have
\[
 \lim_{R\to\infty} \frac1{R^2}\, \int_{R\bD} \bigl( \log_+ \log_+ |F(z)| \bigr)^{1+\e}\,{\rm d}A(z)
= \bE\bigl[ (\log_+\log_+ |F(0)|)^{1+\e} \bigr]\,,
\]
and since the measure $\la$ does not charge constant functions,
by Lemma~\ref{lemma-loglog}, the limit on the RHS is infinite. \hfill $\Box$

\subsection{Proof of Lemma~\ref{lemma-loglog}}

We let $M_u(R)=\max_{R\,\bar\bD} u$ and choose $N$ so that $b^N < R \le b^{N+1}$,
with some $b>1$ to be chosen.
For $1\le j \le N$, we take $z_j$, $|z_j|=R_j$, so that
\[
u(z_j)=M_u(R_j)=b^j, \quad j\in\bN\,,
\]
and let $R_{N+1}=R$. Then, we put
$ \rho_j= R_{j+1}-R_j$, $1\le j \le N$, let $D_j$ be the
disks centered at $z_j$ of radius $\tfrac12 \rho_j$,
and let $ D_j^+ = D_j \bigcap \{R_j \le |z| \le R_{j+1} \}$.
Note that the sets $D_j^+$ are disjoint and that
$A(D_j^+) \ge \tfrac12 A(D_j)$. We claim that
\begin{itemize}
\item {\em If $b$ is chosen sufficiently close to $1$ and $c$ is sufficiently small, then for every $1\le j \le N$, $u(z)\ge c u(z_j)$ on a subset $D_j'\subset D_j^+$ with
$A(D_j') \ge \tfrac14 A(D_j)$}.
\end{itemize}
Indeed, let $D_j^*=\{z\in D_j\colon u(z)\ge cu(z_j)\}$. Then,
\begin{align*}
b^j = u(z_j) &\le \frac1{A(D_j)}\, \int_{D_j} u\, {\rm d}A \\ \\
&\le \frac1{A(D_j)}\, \bigl( c b^j (A(D_j)-A(D_j^*)) + b^{j+1}A(D_j^*) \bigr) \\ \\
&= (b^{j+1}-cb^j)\,\frac{A(D_j^*)}{A(D_j)} + cb^{j}\,,
\end{align*}
whence
\[
\frac{A(D_j^*)}{A(D_j)} \ge \frac{1-c}{b-c} > \frac34
\]
provided that $b>c$  and $3b+c<4$. It remains to put $D_j' = D_j^* \bigcap D_j^+$.

\medskip Now,
\begin{align*}
\int_{R\bD} (\log_+ u)^{1+\e}\, {\rm d}A
&\ge \sum_{j=1}^N \int_{D_j'} (\log_+ u)^{1+\e}\, {\rm d}A \\ \\
&\gtrsim \sum_{j=1}^N j^{1+\e} A(D_j') \\ \\
&\gtrsim \sum_{j=1}^N j^{1+\e} \rho_j^2\,.
\end{align*}
Let $N_0 \gg 1$. Then, for $N \gg N_0$,
\begin{align*}
\Bigl( \sum_{j=N_0}^N \rho_j \Bigr)^2
&< \sum_{j\ge N_0} \frac1{j^{1+\e}} \cdot \sum_{j=N_0}^N j^{1+\e}\rho_j^2 \\ \\
&\lesssim N_0^{-\e}\, \sum_{j=N_0}^N j^{1+\e}\rho_j^2\,.
\end{align*}
Letting $N\to\infty$, we get
\[
\liminf_{R\to\infty} \frac1{A(R\bD)} \int_{R\bD} (\log_+ u)^{1+\e}\, {\rm d}A
\gtrsim N_0^\e\,.
\]
Then, letting $N_0\to\infty$, we conclude the proof. \hfill $\Box$

\section{Proof of Theorems~4 and~5}
Both proofs are quite straightforward.

\subsection{Proof of Theorem~4}

As in the previous proofs we may assume that the measure $\la$ is ergodic.
Then, by the pointwise ergodic theorem, for $\la$-a.e. meromorphic function
$F$,
\[
\lim_{R\to\infty} \frac1{A(R\,\bD)}\, \int_{R\bD} F^{\#}(z)^2\, {\rm d}A(z)
= \bE\bigl[ F^{\#}(0)^2 \bigr]\,.
\]
Since $\la$-a.s., the function $F$ is not a constant (and the distribution of
$\la$ is translation invariant), the RHS is positive
(may be infinite). Thus, for sufficiently large $R$s,
\[
\int_{R\bD} F^{\#}(z)^2\, {\rm d} A(z) \gtrsim R^2\,,
\]
and therefore, $T_F(R)\gtrsim R^2$. \hfill $\Box$

\subsection{Proof of Theorem~5}

Let $F$ be a non-constant locally uniformly recurrent meromorphic function.
We fix a disk $D$ such that $F$ {\em is analytic} on $\bar D$, take
the closed spherical disk $\bar{\mathfrak  D}\subset F(D)$ such that
$\bar{\mathfrak D} \cap F(\partial D)=\emptyset$, and denote by
$\delta$ the spherical distance between $\bar{\mathfrak D}$ and the curve
$F(\partial D)$.

By the definition of local uniform recurrency,
each square $Q$ with sufficiently large length-side $L(Q)$ contains a point $w$
such that
$ \max_{\bar D} \rho (F, \tau_w F) < \tfrac12\, \delta $,
where $\rho$ is the spherical metric.
Denote by $ D_w $ the disk centered at $w$ of the same radius as $D$.
We claim that $\mathfrak D \subset F(D_w)$. To show this, fix a point
$\zeta\in\mathfrak D\subset F(D)$. Then, by the argument principle,
the index of the curve $F(\partial D)$ with respect to the point $\zeta$ is positive.
Furthermore, when the point $z$ traverses the circumference $\partial D$,
$F(z)$ traverses the curve $F(\partial D)$, $F(z+w)$ traverses the curve
$F(\partial D_w)$, and the spherical
distance between $F(z)$ and $F(z+w)$ remains less than $\frac12\,\delta $,
while $ \rho(\zeta, F(\partial D))\ge \delta $. Hence, the index of the curve
$F(\partial D_w)$ with respect to $\zeta$ coincides with that of $F(\partial D)$,
and therefore, is positive as well. Thus, $\zeta\in F(D_w)$ proving the claim.

Denote by $Q^*$ the square having the
same center as $Q$ and with the length-side $L(Q^*)=L(Q)+\operatorname{radius}(D)$.
Since $D_w\subset Q^*$, we conclude that $\mathfrak D\subset F(Q^*)$. Hence,
the spherical area of $F(Q^*)$ is not less than that of $\mathfrak  D$.
Packing the disk $R \bD=\{|z|<R\}$  with sufficiently large $R$ by about $cR^2$ disjoint translations of the square $Q^*$, we see that the
spherical area of $F(R\bD)$ is $\gtrsim R^2$, which yields the theorem.
\mbox{} \hfill $\Box$

\section{Entire functions of almost minimal growth outside a ternary system of squares}

\subsection{Ternary system of squares}

We will construct the closed set $E\subset \bC$ which we will call
{\em the ternary systems of squares}. It will be defined as the limit
of the increasing sequence $(E_n)$ of compact sets such that $E_n$ consists
of $E_{n-1}$ and its eight disjoint translations. One can think about the
limiting set $E$ as a two-dimensional ternary Cantor-type set viewed from  the
inside-out.

\subsubsection{Notation}

For $X\subset\bC$ and $\eta>0$, we put
\[
X^{+\eta}=\bigl\{z\in\bC\colon d_\infty (z, X) \le \eta \bigr\}, \quad
X^{-\eta}=\bigl\{z\in\bC\colon d_\infty (z, X^c) \ge \eta \bigr\}\,.
\]
Here and elsewhere, $d_\infty$ denotes the $\ell^\infty$-distance on $\bR^2$.

\medskip
For $X\subset\bC$, we put
$\tau_w X = \bigl\{z-w\colon z\in X \bigr\}$. That is, if the function $f$
is defined on $X$, then $\tau_w f$ is defined on $\tau_w X$.

\subsubsection{Squares and corridors}
\begin{figure}[h!]
\centering
\includegraphics[scale=0.75]{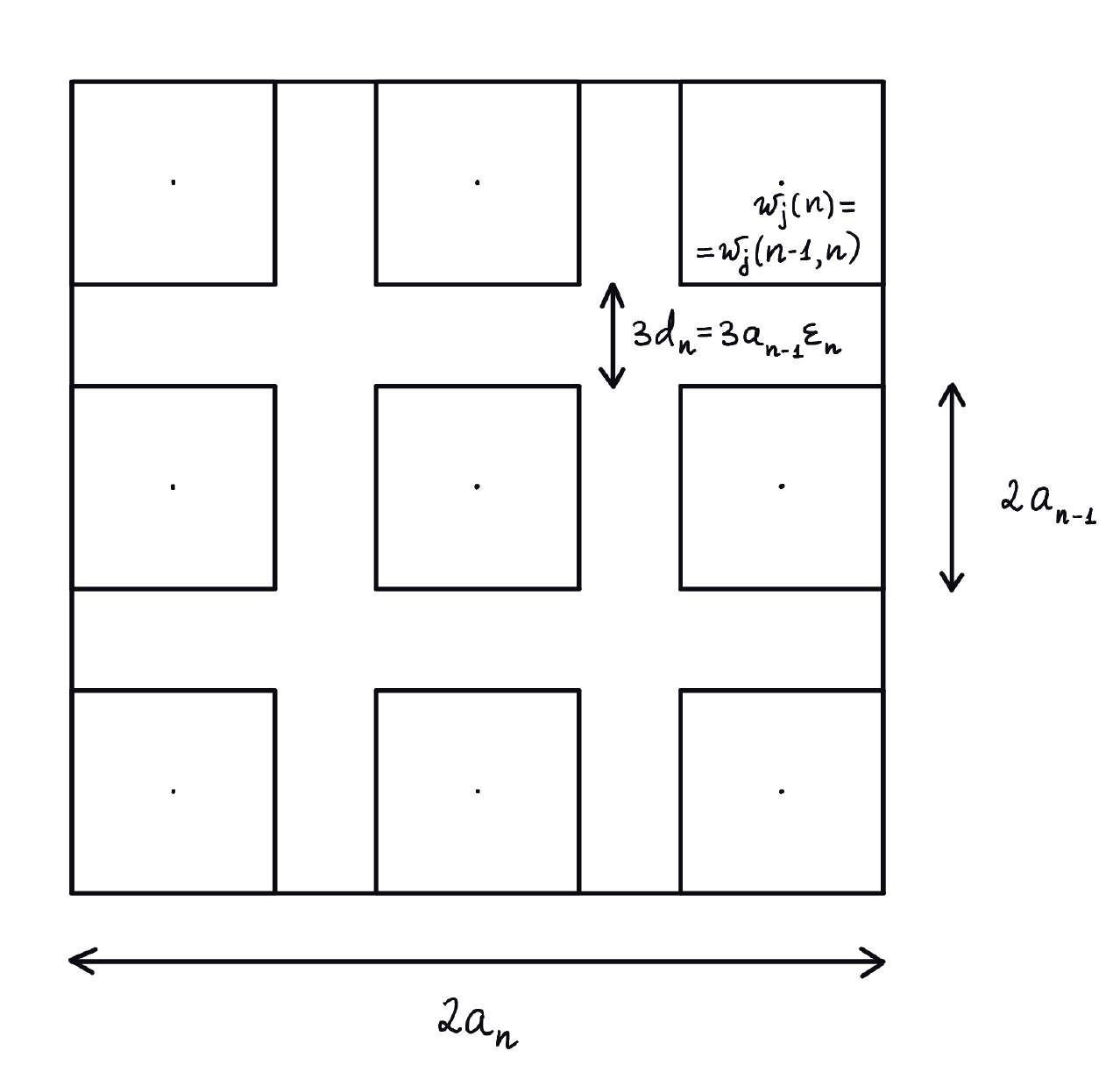}
\caption{\label{fig:squares}}
The square $S_n$ with 9 copies of the square $S_{n-1}$.
\end{figure}

We fix a sequence $(\e_n)\downarrow 0$ and define:
\begin{itemize}
\item the increasing sequence $(a_n)$ by $a_0=1$, $a_n=3a_{n-1}(1+\e_n)$;
\item the squares $ S_n=[-a_n, a_n]^2$;
\item the translates $w_j(n)=a_{n-1}(2+3\e_n)\omega_j$, where
$\{\omega_j\}_{0\le j \le 8} =\{0, \pm 1, \pm {\rm i}, \pm 1\pm {\rm i}\}$;
\end{itemize}
see Figure~\ref{fig:squares}.
Then, $E_0=S_0$, $ E_n=\bigcup_{j=0}^8 \tau_{w_j(n)}E_{n-1} $, and finally,
$E=\bigcup_n E_n$.

\begin{figure}[h!]
\centering
\includegraphics[scale=0.75]{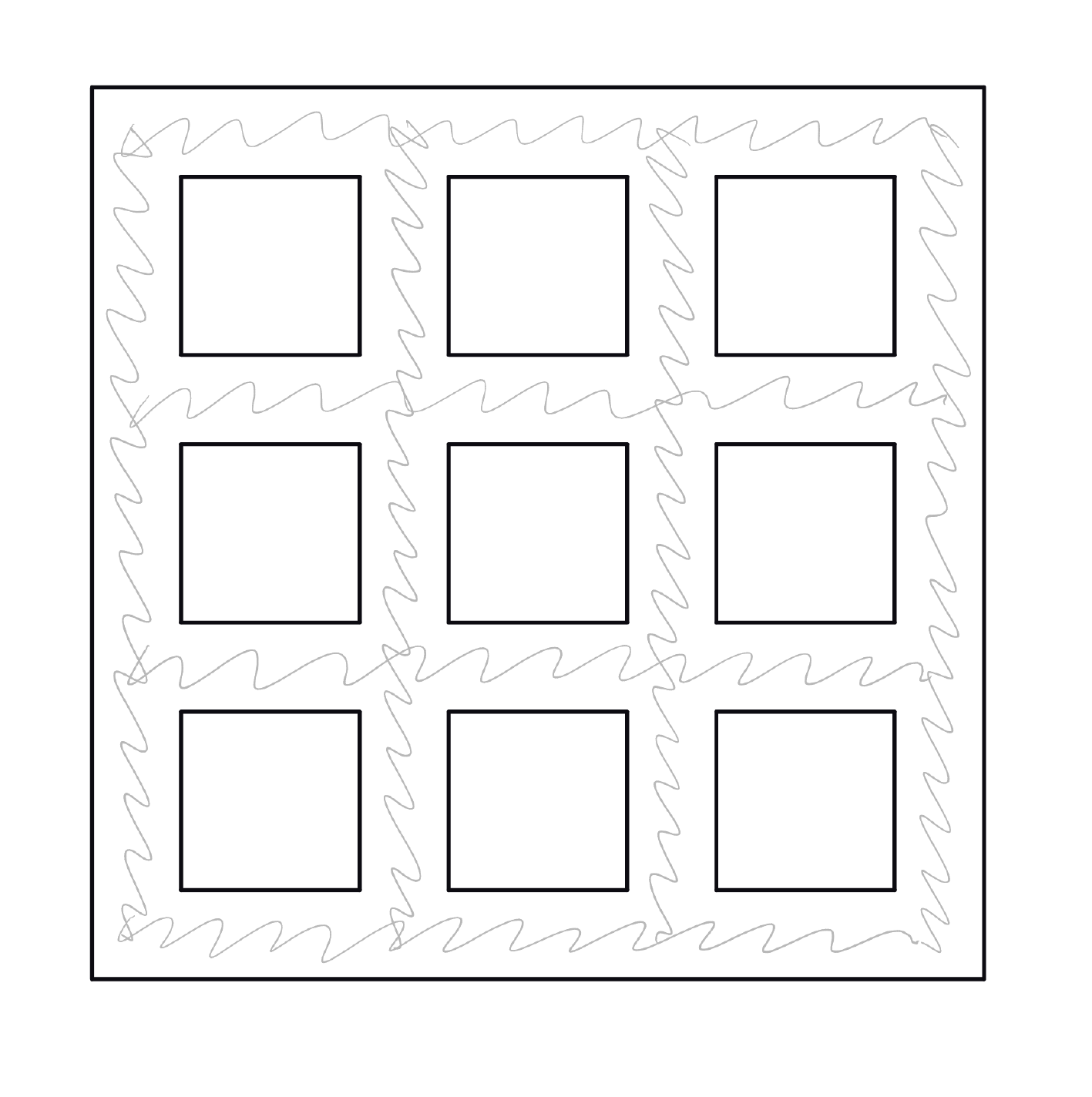}
\caption{\label{fig:corridors}}
The corridors $K_n$.
\end{figure}

\medskip
For every $k<n$, the set $E_n$ consists of $9^{n-k}$ disjoint copies of $E_k$. We denote
by $w_j(k, n)$, $0\le j \le 9^{n-k}-1$, the centers of these copies. That is,
there exist indices $j_1, \ldots , j_{n-k}\in \{0, 1, \ldots , 8\}$ such that
\[
w_j(k, n)= w_{j_1}(k+1)+ \ldots +w_{j_{n-k}}(n)
\]
(in particular, $w_j(n-1,n)=w_j(n)$), and
\[
E_n = \bigcup_{j_n=0}^8 \ldots \bigcup_{j_1=0}^8 \tau_{w_{j_n}(n)+\ldots+w_{j_1}(1)} E_0
= \bigcup_{j=0}^{9^{n-k}-1} \tau_{w_j(k, n)} E_{k}\,.
\]
Next, we denote by $K_n$
the union of the corridors left on the $n$th step of the construction and the outer perimeter
corridor that goes along the boundary $\partial S_n$ (see Figure~\ref{fig:corridors}).
The width of these corridors is $3d_n$, where $d_n=a_{n-1}\e_n$. That is,
\[
K_n = S_n^{+3d_n} \setminus \bigcup_{j=0}^8 \tau_{w_j(n)}S_{n-1}\,.
\]

To simplify computations,
in what follows, we always assume that $\e_1  <1$ and that
$\e_n \ge \e_{n+1}\ge \frac13 \e_n$. Since
\[
\frac{d_{n+1}}{d_n} = \frac{3(1+\e_n)\e_{n+1}}{\e_n}\,,
\]
these assumptions yield that
\[
1 < \frac{d_{n+1}}{d_n} < 6\,.
\]
In particular, the sequence $(d_n)$ is increasing.

\subsubsection{Fat systems of squares}

We will call the set $E$ {\em fat} if $\sum_{n\ge 1} \e_n<\infty$.
In this case,
\[
a_n = 3^n \prod_{j=1}^n (1+\e_j) = (a+o(1)) 3^n
\]
with $a>0$.

Note that since $E_n$ consists of $9^n$ disjoint translations of the square $S_0$,
$A(E_n) = 4 \cdot 9^n$. If the set $E$ is fat then $A(S_n)=(a^2+o(1))9^n$, and
$A(E_n)/A(S_n)=b+o(1)$ with some $b>0$. In particular, fat sets $E$ (and only they)
have positive relative area.

\subsubsection{Concordance and $\delta$-concordane}

Given a ternary system of squares $E$, we call the function
$\Phi\colon \bC\to\bC$ {\em concordant} with $E$ if
for every $n>k\ge 1$ and $0\le j \le 9^{n-k}-1$,
\[
\tau_{w_j(k, n)}\Phi = \Phi \quad {\rm everywhere\ on\ } S_k.
\]
Given a sequence
$\delta=(\delta_k)\downarrow 0$, we say that the function
$\Phi$ is $\delta$-concordant with $E$ if
for every $n>k\ge 1$ and $0\le j \le 9^{n-k}-1$,
\[
\max_{S_k} \bigl| \tau_{w_j(k, n)} \Phi - \Phi \bigr| < \delta_k\,.
\]

\subsection{Main Lemma}

For a continuous function $\Phi$ and a compact set $K$, we put
\[
M_\Phi (K) = \max_{K} |\Phi|, \quad  m_{\Phi}(K) = \min_K |\Phi|\,.
\]

Define {\em the majorant}
\[
\mathcal M_B(n) = \exp\Bigl( Bn + \pi \sum_{j=1}^{n} \tfrac1{\e_j} \Bigr), \quad n\ge 1
\]
with sufficiently large positive $B$ and put $\mathcal M_B(0)=1$. Then, define the sequence $\Delta$ by
\[
\Delta_n  = \exp\Bigl( -\frac1{10} \mathcal M_B(n-1)\Bigr), \qquad n\ge 1\,.
\]
\begin{lemma}\label{lemma-EF}
For any sufficiently large positive $B$,
there exists a non-constant entire function $G$ which is $\Delta$-concordant with $E$, satisfies
\[
\log M_G\bigl( S_n \bigr) \lesssim e^{-B}
\mathcal M_B(n)\,,
\]
and
\[
\max_{S_0} |G(z)-z| \le \frac13\,.
\]
\end{lemma}

We start with the subharmonic counterpart of this lemma.

\subsection{Subharmonic construction}

\begin{lemma}\label{lemma-SH}
For any $B\ge  20$, there exists a sequence of non-negative subharmonic functions $u_n$ in $\bC$ with the following properties:

\smallskip\noindent{\rm (i)} for each $j\in\{0, 1, \ldots , 8\}$,
\[
\tau_{w_j(n-1, n)}u_{n-1} = u_n \quad \text{on}\  S_{n-1};
\]

\smallskip\noindent{\rm (ii)} $M_{u_n}(S_n^{+d_{n+1}}) < e^{-B+10}\mathcal M_B(n)$;

\smallskip\noindent{\rm (iii)} $m_{u_n}(K_n^{-\frac12 d_n}) > \tfrac12\, \mathcal M_B(n-1)$.
\end{lemma}

Note that by property (i), for every $m>n$, $u_m=u_n$ everywhere on
the square $S_n$. Hence, the sequence $(u_n)$ converges to the limiting subharmonic
function $u$. By property (i), the limiting function is concordant with $E$.
By (ii), we have $M_u(S_n) \le e^{-B+10} \mathcal M_B(n)$, $n\ge 0$, and by (iii),
$m_{u}(K_n^{-\frac12 d_n}\bigcap S_n) \le \tfrac12\, \mathcal M_B(n-1)$, $n\ge 1$.

\medskip
\noindent{\em Proof of Lemma~\ref{lemma-SH}}:
Take the subharmonic function
\[
h(z) =
\begin{cases}
\cosh x \cos y\, \ & |y|<\tfrac{\pi}2, \\
0 \ &\text{otherwise}\,,
\end{cases}
\]
scale it
\[
h_n(z) = h\Bigl( \frac{\pi}{3d_n} z\Bigr)\,,
\]
\[
\xi_n = a_{n-1} + \frac32 d_n = a_{n-1}\Bigl( 1 + \frac32 \e_n \Bigr)\,,
\]
\begin{figure}[h!]
\centering
\includegraphics[scale=0.75]{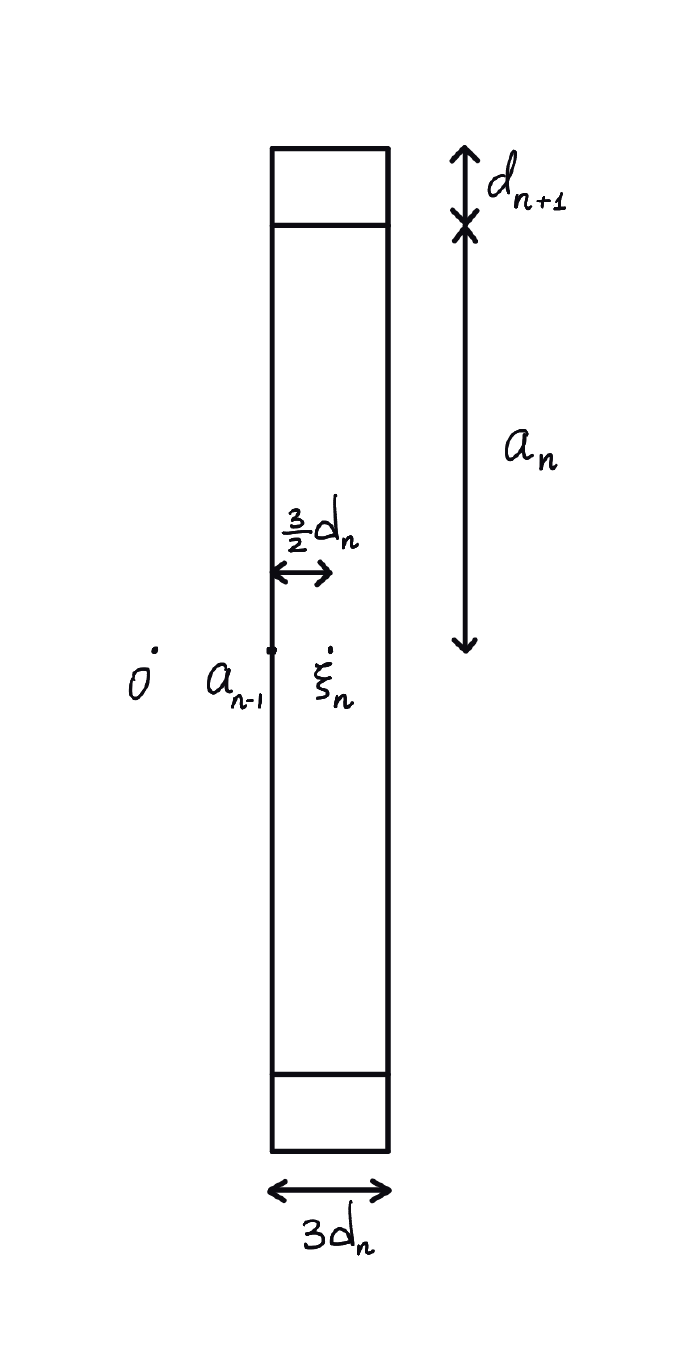}
\caption{\label{fig:h}}
Scaling, shifting and rotating the function $h$.
\end{figure}
and take the upper envelope of $8$ shifted and rotated copies of $h_n$:
\begin{align}\label{eq_v}
\nonumber
v_n(z)=\max\bigl\{ &h_n(z+{\rm i}\xi_n),  h_n(z-{\rm i}\xi_n),
h_n({\rm i}(z+\xi_n)), h_n({\rm i}(z-\xi_n)), \\
&h_n(z+3{\rm i}\xi_n),  h_n(z-3{\rm i}\xi_n),
h_n({\rm i}(z+3\xi_n)), h_n({\rm i}(z-3\xi_n))
\bigr\}\,.
\end{align}
We will need two estimates:
\begin{equation}\label{eq-LB-v_n}
m_{v_n} \bigl( K_n^{-\frac12 d_n} \bigr) =
\cos\Bigl( \frac{\pi}{3d_n} \cdot d_n \Bigr) = \frac12\,,
\end{equation}
and
\begin{align}\label{eq-UB-v_n}
\nonumber M_{v_n} \bigl( S_n^{+d_{n+1}} \bigr)
&\le \exp\Bigl( \frac{\pi}{3d_n}\, \bigl( a_n + d_{n+1} \bigr) \Bigr)
\stackrel{d_{n+1}<6d_n}< \exp\Bigl( \frac{\pi}{\e_n}\, \frac{a_n}{3a_{n-1}} + 2\pi  \Bigr) \\ \nonumber \\
&= \exp\Bigl( \frac{\pi}{\e_n}\, (1+\e_n) + 2\pi \Bigr)
= \exp\Bigl( \frac{\pi}{\e_n}\, + 3\pi \Bigr)\,.
\end{align}

We put $u_0=\tfrac13$. Assuming that the subharmonic functions $u_0, \ldots , u_{n-1}$ have
been already defined, we glue together the functions $\tau_{w_j(n-1, n)} u_{n-1}$, putting
\[
u_n =
\begin{cases}
\max\bigl\{ \mathcal M_B(n-1) v_n, \tau_{w_j(n-1, n)} u_{n-1} \bigr\}
\ &\text{on\ } \tau_{w_j(n-1, n)} S_{n-1}^{+\frac12 d_n}, \ 0\le j \le 8\,, \\
\mathcal M_B(n-1) v_n \ &\text{otherwise}\,,
\end{cases}
\]
where $v_n$ is the subharmonic function defined in~\eqref{eq_v}.
Note that this definition ensures property (i) in the statement of the lemma, as
\[
v_n= 0 \quad \text{on}\ \bigcup_{j=0}^8 \tau_{w_j(n-1, n)} S_{n-1}\,.
\]

We claim that, for $B \ge 20$ and $n\ge 1$,
\[
\max_{\partial (S_{n-1}^{+\frac12 d_n})} u_{n-1} <
\mathcal M_B(n-1) \min_{\partial (S_{n-1}^{+\frac12 d_n})} v_n
\]
(with $\mathcal M_B(0)=1$).
This claim yields that
\[
\tau_{w_j(n-1, n)} u_{n-1} < \mathcal M_B(n-1) v_n \quad {\rm on\ }
\partial\bigl( \tau_{w_j(n-1, n)} S_{n-1}^{+\frac12 d_n} \bigr),
\]
and therefore, the functions $u_n$, $n\ge 1$, are subharmonic in $\bC$.

\medskip
The case $n=1$ of our claim follows from the lower bound
$v_1\ge \tfrac12 $ on $\partial (S_0^{+\frac12 d_1})$.
Now, let $n\ge 2$. We know that $ u_{n-1} = \mathcal M_B(n-2) v_{n-1} $
outside the set $ \bigcup_{j=0}^8 \tau_{w_j(n-2, n-1)} S_{n-2}^{+\frac12 d_{n-1}} $.
Note that
\[
\bigcup_{j=0}^8 \tau_{w_j(n-2, n-1)} S_{n-2}^{+\frac12 d_{n-1}}
\ \subset \ S_{n-1}^{+\frac12 d_{n-1}}
\ \stackrel{d_{n-1}<d_n}\subset \ {\rm interior} \bigl(S_{n-1}^{+\frac12 d_n}\bigr).
\]
Hence,
$ u_{n-1} = \mathcal M_B(n-2) v_{n-1} $ on $\partial \bigl(S_{n-1}^{+\frac12 d_n}\bigr)$.
Furthermore, by the bound~\eqref{eq-UB-v_n}, on
$\partial \bigl(S_{n-1}^{+\frac12 d_n}\bigr)$ we have
$ \mathcal M_B(n-2) v_{n-1} < \mathcal M_B(n-2) \cdot e^{\pi/\e_{n-1} +3\pi}
= e^{-B+3\pi} \mathcal M_B(n-1) $.
For $B \ge 20> 3\pi + \log 2$, we have $e^{-B+3\pi} < \tfrac12$, whence
$ u_{n-1} < \frac12 \mathcal M_B(n-1) $ on $\partial \bigl(S_{n-1}^{+\frac12 d_n}\bigr)$.
On the other hand,
$ \partial \bigl(S_{n-1}^{+\frac12 d_n}\bigr) \subset K_n^{-\frac12 d_n } $,
so applying the lower bound~\eqref{eq-LB-v_n} for $v_n$, we get the claim.

\medskip
Note that $u_n = \mathcal M_B(n-1) v_n$ on
$\partial \bigl( S_n^{+d_{n+1}} \bigr)$, and
by~\eqref{eq-UB-v_n}, \[ \mathcal M_B(n-1) v_n < e^{-B+3\pi} \mathcal M_B(n) \]
therein. This proves (ii). At last, on $ K_n^{-\frac12 d_n}$, we have
\[ u_n=\mathcal M_B(n-1) v_n > \tfrac12\, \mathcal M_B(n-1)\,, \] proving (iii).
\hfill $\Box$

\subsection{Proof of Lemma~\ref{lemma-EF}}

\subsubsection{Beginning the proof}
We put $G_1(z)=z$ and
construct a sequence $(G_n)$ of entire functions with the following properties:

\smallskip\noindent{\rm (i)} for $n\ge 2$ and $j\in \{0, 1, \ldots , 8\}$,
\[
\max_{S_{n-1}} \bigl| G_{n-1} - \tau_{w_j(n-1, n)} G_n \bigr| < \frac1{10}\, \Delta_n\,,
\]

\smallskip\noindent{\rm (ii)} for $n\ge 2$,
\[
\log M_{G_n}\bigl( S_n^{+\frac9{10}\, d_{n+1}} \bigr) <
e^{-B+10} \mathcal M_B(n)\,.
\]
Then, the existence of the function $G$ will follow from the following claim.

\begin{claim}\label{claim2}
For every $1\le k <n$,
\[
\max_{S_k} \bigl|  G_k - \tau_{w_j(k, n)} G_n \bigr|
< \frac1{10}\, \sum_{i=k+1}^n  \Delta_i\,.
\]
\end{claim}
First, assuming that estimates (i) and (ii) and the claim hold, we complete the proof of Lemma~\ref{lemma-EF}. On the second step, we prove the claim assuming that the property
(i) holds. On the last step, we construct the sequence $(G_n)$ having properties
(i) and (ii).

\medskip
We put
\[
G=G_1+\sum_{i\ge 2} ( G_{i} - G_{i-1} )\,.
\]
By (i), the series converges locally uniformly in $\bC$. Moreover,
\[
\max_{S_k} |G_k-G| \le \frac1{10}\, \sum_{i=k+1} \Delta_i
\]
and then, for $n\ge k$,
\[
\max_{S_k} |G_n-G| \le
\max_{S_n} |G_n-G| \le
\frac1{10}\, \sum_{i=n+1} \Delta_i\,.
\]
Combining these inequalities with the claim, we conclude that, for every $n> k \ge 1$,
\begin{align*}
\max_{S_k} |G-\tau_{w_j(k, n)}G|  &\le  \Bigl( \frac1{10}\, \sum_{i=k+1} \Delta_i \Bigr)
+ \Bigl(  \frac1{10}\, \sum_{i=n+1} \Delta_i \Bigr)
+\Bigl( \frac1{10}\, \sum_{i=k+1} \Delta_i \Bigr) \\
&< \sum_{j=k+1} \Delta_i < \Delta_k
\end{align*}
provided that the parameter $B$ is large enough. That is, the limiting entire function
$G$ is $\Delta$-concordant with $E$.

Furthermore, by properties (i) and (ii), the function $G$ satisfies
$ \log M_G(S_n) \lesssim e^{-B} \mathcal M_B(n) $ and
\[
\max_{S_0} |G(z)-z| \le \max_{S_1} |G(z)-G_1(z)| \le \frac1{10}\, \sum_{i\ge 2} \Delta_i < \frac13\,,
\]
provided that $B$ is sufficiently large.

\subsubsection{Proof of Claim~\ref{claim2}}
We use induction on $n-k$.  The base of induction $n-k=1$ is exactly  property (i).

Now, assuming that the claim holds for the pair $(k, n-1)$, we will prove it for $(k, n)$.
For every $0\le j \le 9^{n-k}$, we write
\[
j=\sum_{\ell=1}^{n-k} j_\ell\, 9^{\ell-1}\,,
\]
and put
\[
j' = \sum_{\ell=1}^{n-k-1} j_\ell\, 9^{\ell-1} = j - j_{n-k}\,9^{n-k-1}\,.
\]
Then we have
\[
\max_{S_k} | G_k - \tau_{w_{j'}(k, n-1)} G_{n-1} |  \le \frac1{10}\, \sum_{i=k+1}^{n-1} \Delta_i
\quad (\text{induction\ hypothesis}) \eqno ({\rm a})
\]
and
\[
\max_{S_{n-1}} | G_{n-1} - \tau_{w_{j_{n-k}}(n-1, n)} G_n |  \le \frac1{10}\, \Delta_n
\quad (\text{property\ (i)}) \eqno ({\rm b})
\]
Then, taking into account that $\tau_{-w_{j'}(k, n-1)}S_k \subset S_{n-1}$ and using (b), we get
\[
\max_{S_k} | \tau_{w_{j'} (k, n-1)} G_{n-1} - \tau_{w_{j}(k, n)} G_n |  \, =
\, \max_{ \tau_{-w_{j'} (k, n-1)} S_k} |  G_{n-1} - \tau_{w_{j_{n-k}}(n-1, n)} G_n |
\]
\[
\le \max_{ S_{n-1}} |  G_{n-1} - \tau_{w_{j_{n-k}}(n-1, n)} G_n  | \le \frac1{10}\, \Delta_n\,.
\eqno ({\rm c})
\]
Now, adding (a) and (c), we conclude proof of the claim. \hfill $\Box$

\medskip
Thus, it remains to construct a sequence of entire functions $(G_n)$ satisfying
conditions (i) and (ii).

\subsubsection{Constructing the sequence $(G_n)$}

We fix a sequence of smooth cut-off functions $\chi_n$, $0\le\chi_n\le 1$, so that
\[
\chi_n = \begin{cases}
1 \ &\text{on\ } S_{n-1}^{+\frac35\, d_n} \\
0 \ &\text{on\ } \bC\setminus S_{n-1}^{+\frac45\, d_n}\,.
\end{cases}
\]
and
$ \sup_n \| \nabla \chi_n \|_\infty < \infty $ (such a sequence exists since
$d_n \ge d_1 > 0$).

We put $G_1(z)=z$ and suppose that the functions $G_1$, \ldots , $G_{n-1}$ have already been
constructed. We put
\[
g_n = \sum_{j=0}^8 \tau_{w_j(n-1, n)} \bigl( \chi_n G_{n-1} \bigr)\,,
\]
and note that
\[
\bar\partial g_n = \sum_{j=0}^8 \tau_{w_j(n-1, n)} \beta_n\,,
\]
with $\beta_n = G_{n-1}\bar\partial \chi_n$.
Then, we define the function $G_n$ by $G_n = g_n -\alpha_n$, where
$\alpha_n$ is H\"ormander's solution~\cite[Theorem~4.2.1]{Hormander}
to the $\bar\partial$-equation $\bar\partial \alpha_n = \bar\partial g_n$ satisfying
\begin{equation}\label{eq-Horm}
\int_{\bC} |\alpha_n|^2 e^{-u_n}\, \frac{{\rm d} A}{(1+|z|^2)^2}
< \frac12\, \int_{\bC} |\bar\partial g_n|^2 e^{-u_n}\, {\rm d}A\,,
\end{equation}
where $u_n$ are the subharmonic functions constructed in Lemma~\ref{lemma-SH}.

\subsubsection{Estimating the integral on the RHS of~\eqref{eq-Horm}}

Note that
\begin{align*}
\operatorname{spt}(\bar\partial g_n) &=
\bigcup_{j=0}^8 \tau_{w_j(n-1, n)} \operatorname{spt}(\bar\partial \chi_n) \\
&\subset \bigcup_{j=0}^8 \tau_{w_j(n-1, n)} \bigl( S_{n-1}^{+\frac45\, d_n }
\setminus\text{interior}(S_{n-1}^{+\frac35\, d_n })  \bigr) \subset K_n^{-\frac12 d_n}\,,
\end{align*}
and therefore, $u_n>\tfrac12\, \mathcal M_B(n-1)$ on $ \operatorname{spt}(\bar\partial g_n) $.
Furthermore, since $ \operatorname{spt}(\bar\partial \chi_n) \subset S_{n-1}^{+\frac45 d_n}$,
we have
\[
| \bar\partial g_n | \le C_\chi M_{G_{n-1}}(S_{n-1}^{+\frac{9}{10} d_n})
< C_\chi \exp\Bigl( e^{-B+10} \mathcal M_B(n-1)\Bigr)
\]
with $C_\chi = \sup_n \| \nabla\chi_n \|_\infty$ (in the second inequality we have
used the inductive assumption). Taking into account that the area of
$\operatorname{spt}(\bar\partial g_n)$ is less than $(a_n+d_n)^2 < (6^n+6^n)^2$
and recalling that $u_n \ge \frac12 \mathcal M_B(n-1)$ on $K^{-\frac12 d_n}$,
we conclude that the integral on the RHS of~\eqref{eq-Horm} does not exceed
\[
4\cdot 6^{2n} C_\chi^2 \exp\Bigl( \bigl( 2e^{-B+10}-\frac12 \bigr) \mathcal M_B(n-1)\Bigr)
<\exp\Bigl( -\frac25 \mathcal M_B(n-1) \Bigr)\,,
\]
provided that the constant $B$ is sufficiently large.

Therefore, by H\"ormander's theorem,
\begin{equation}\label{eq-norm-alpha}
\int_{\bC} |\alpha_n|^2 e^{-u_n}\, \frac{{\rm d} A}{(1+|z|^2)^2}
<\frac12\, \exp\Bigl( -\frac25 \mathcal M_B(n-1) \Bigr)\,.
\end{equation}

\subsubsection{Proving property (ii) for the sequence $G_n$}
Here, we aim to show that
\[
\log M_{G_n} \bigl( S_n^{+\frac9{10} d_{n+1}} \bigr) < e^{-B+10}\, \mathcal M_B(n)\,.
\]
Let $ c_0<\frac1{10}d_1 $ be a positive constant.
Then, for $z\in S_n^{+\frac9{10}d_{n+1}}$,
we have
\[
|G_n(z)|^2 \le \frac1{\pi c_0^2}\, \int_{\tau_z (c_0\bD)} |G_n|^2
\le \frac2{\pi c_0^2}\, \int_{\tau_z (c_0\bD)} (|g_n|^2 + |\alpha_n|^2)\,.
\]
To estimate the first integral, we observe that
\[
\| g_n \|_\infty \le \max_{S_{n-1}^{+\frac45\, d_n}} |G_{n-1}|
< \exp\Bigl( e^{-B+10} \mathcal M_B(n-1)\Bigr)
\]
(in the second inequality we have used the induction assumption). Thus,
\[
\frac2{\pi c_0^2}\, \int_{\tau_z (c_0\bD)} |g_n|^2
< \frac12\, \exp\Bigl( 2e^{-B+10} \mathcal M_B(n) \Bigr)\,,
\]
provided that the constant $B$ is sufficiently large.

To estimate the second integral, using the fact
that $z\in S_n^{+\frac9{10}d_{n+1}}$, we write
\begin{align*}
\int_{\tau_z (c_0\bD)} |\alpha_n|^2
&< \int_{\bC} |\alpha_n|^2 e^{-u_n}\, \frac{{\rm d} A}{(1+|z|^2)^2}
\cdot C(a_n+d_{n+1})^4\, \exp\Bigl( \max_{S_n^{+d_{n+1}}} u_n \Bigr) \\ \\
&< C_1 6^{4n}\, \exp\Bigl( -\frac25 \mathcal M_B(n-1) + e^{-B+20} \mathcal M_B(n) \Bigr)\,,
\end{align*}
whence,
\[
\frac2{\pi c_0^2}\, \int_{\tau_z (c_0\bD)} |\alpha_n|^2
< \frac12 \exp\Bigl( 2e^{-B+10} \mathcal M_B(n) \Bigr)\,,
\]
again, provided that the constant $B$ is large enough. Thus,
\[
|G_n| < \exp\Bigl( e^{-B+10} \mathcal M_B(n) \Bigr)
\]
everywhere on $S_n^{+\frac9{10}d_{n+1}}$, as we have claimed.

\subsubsection{Proving property (i) for the sequence $G_n$}
First, we note that
\[
\max_{S_{n-1}} | G_{n-1} - \tau_{-w_j(n-1, n)} G_n|
= \max_{\tau_{w_j(n-1, n)} S_{n-1}} |\alpha_n|\,,
\]
and that $\alpha_n$ is analytic in the $c_0$-neighbourhood of
each of the sets $ \tau_{w_j(n-1, n)} S_{n-1} $. Then,
for every $j\in\{0, 1, \ldots , 8\}$ and
every $z\in\tau_{w_j(n-1, n)} S_{n-1}$, we have
\begin{align*}
|\alpha_n(z)|^2 &\le \frac1{\pi c_0^2} \int_{\tau_z(c_0\bD)} |\alpha_n|^2 \\ \\
&< \int_{\bC} |\alpha_n|^2 e^{-u_n}\, \frac{{\rm d} A}{(1+|z|^2)^2}
\cdot C(a_n+c_0)^4 \exp\Bigl( \max_{\tau_{w_j(n-1, n)} S_{n-1}^{+c_0}}  u_n \Bigr) \\ \\
&\stackrel{\eqref{eq-norm-alpha}}< C' (a_n +c_0)^4
\exp\Bigl( \, \max_{\tau_{w_j(n-1, n)} S_{n-1}^{+c_0}}\, u_n
- \frac25 \mathcal M_B(n-1) \Bigr)\,.
\end{align*}
Recall that $u_n = \tau_{w_j(n-1, n)} u_{n-1}$ on
each square $ \tau_{w_j(n-1, n)} S_{n-1} $. Then, recalling that $d_n>1$ and choosing $c_0<1$, we see
that
\[
\max_{\tau_{w_j(n-1, n)}S_{n-1}^{+c_0}} u_n = \max_{S_{n-1}^{+c_0}} u_{n-1}
\le  \max_{S_{n-1}^{+d_n}} u_{n-1}
< e^{-B+10} \mathcal M_B(n-1)\,.
\]
Therefore,
\[
|\alpha_n(z)|^2 \le C_1 6^{4n} \exp\Bigl( -\frac25 \mathcal M_B(n-1)
+ e^{-B+10} \mathcal M_B(n-1) \Bigr)
< \frac1{100}\, e^{-\frac15\, \mathcal M_B(n-1)}\,,
\]
provided that $B$ is sufficiently large, and finally,
\[
\max_{\tau_{w_j(n-1, n)} S_{n-1}} |\alpha_n|
< \frac1{10}\, e^{-\frac1{10}\, \mathcal M_B(n-1)} = \frac1{10}\, \Delta_n\,,
\]
again, provided that $B$ is sufficiently large. This completes the (somewhat long) proof
of Lemma~\ref{lemma-EF}. \hfill $\Box$

\section{A version of the Krylov-Bogolyubov construction}

\subsection{Some notation}
In this section, we denote by $(S_n)$ any increasing sequence of squares centered at the origin with the side-lengths tending to infinity.

\medskip
If $S\subset\bC$ is a square and $X\subset\bC$ is a Borel set, then
we denote the relative area of $X$ in $S$ by
\[
A_S(X)=\frac{A(X\bigcap S)}{A(S)}\,.
\]

\medskip
For an entire function $G\in \cE$, let $\mathcal O_G=\{\tau_w G\}_{w\in\bC}$
denote its orbit and $\bar{\mathcal O}_G$ denote the closure of $\mathcal O_G$
in $\cE$.

\medskip
For a compact set $K\subset\bC$ and a continuous function $f\colon K\to \bR$, we denote
by $\operatorname{osc}_K f = \max_K f - \min_K f $, the oscillation of $f$ on $K$.

\subsection{The Lemma}

\begin{lemma}\label{lemma-KB}
Let $G\in\cE$.

\smallskip\noindent{\rm (i)}
Suppose that there exists an increasing sequence $(M_k)\uparrow +\infty$ such that
\begin{equation}\label{eq-KB1}
\lim_{k\to\infty}\, \liminf_{n\to\infty}\,
A_{S_n}\bigl\{w\colon \max_{\tau_w S_k} |G| \le M_k \bigr\} = 1\,,
\end{equation}
and there exists a square $S$ and a constant $c>0$ such that
\begin{equation}\label{eq-KB2}
\limsup_{n\to\infty} A_{S_n}\bigl\{w\colon \osc_{\tau_w S} |G| \ge c \bigr\}>0\,.
\end{equation}
Then there exists a translation-invariant probability measure $\la$ supported by
$\bar{\mathcal O}_G$ which does not charge the constant functions.

\smallskip\noindent{\rm (ii)} Furthermore, suppose that condition~\eqref{eq-KB1}
is replaced by a stronger one:
\begin{equation}\label{eq-KB1'}
\sum_{k\ge 1} \Bigl( 1- \liminf_{n\to\infty}\,
A_{S_n}\bigl\{w\colon \max_{\tau_w S_k} |G| \le M_k \bigr\} \Bigr) < \infty
\end{equation}
and that condition~\eqref{eq-KB2} continues to hold. Then, for $\la$-a.e. $F\in\cE$,
\[
\limsup_{k\to\infty} \bigl( \max_{S_k} |F| - M_k \bigr) \le 0\,.
\]
\end{lemma}

It is worth mentioning that condition (i) already yields the upper bound though only
on a subsequence of the squares $S_k$: for $\la$-a.e. $F\in\cE$,
\[
\liminf_{k\to\infty} \bigl( \max_{S_k} |F| - M_k \bigr) \le 0\,.
\]

\subsection{Proof of part (i) of Lemma~\ref{lemma-KB}}
Consider the sequence of probability measures on $\cE$:
\[
\la_n = \frac1{A(S_n)}\, \int_{S_n} \delta_{\tau_w G}\, {\rm d} A(w)\,.
\]
In other words, for any Borel set $\mathcal X\subset\cE$,
\[
\la_n(\mathcal X) =
\frac1{A(S_n)}\, \int_{S_n} \done_{\mathcal X} (\tau_w G)\, {\rm d} A(w)
=A_{S_n} \bigl\{w\colon \tau_w G \in\mathcal X \bigr\}.
\]

\subsubsection{Tightness of the sequence $(\la_n)$}
We claim that $(\la_n)$ {\em is a tight sequence} of probability measures,
that is, for every $\delta>0$, there exists a compact set $\mathcal K\subset\cE$ such that,
for every $n\ge 1$, $\la_n(\mathcal K)>1-\delta$.

To see this, given $k\ge 2$, we choose $n_k>k$ so that
for $n\ge n_k$,
\[
A_{S_n}\bigl\{w\colon \max_{\tau_w S_k} |G| > M_{n_k} \bigr\} < \frac1{k^2}\,,
\]
and let
\[
\mu_k = \max_{2S_{n_k-1}} |G| + M_{n_k}\,,
\]
where $2S_{n_k-1}$ is the square concentric with $S_{n_k-1}$ and having
double the side-length. The sets
\[
\mathcal K_\ell = \bigl\{F\in\cE\colon \max_{S_k} |F| \le \mu_k \ {\rm for\ } k\ge \ell\bigr\}
\]
are compact subsets of $\cE$. We will show that for any $n\ge 1$ and any $\ell\ge 2$,
\[
\la_n (\mathcal K_\ell) \ge 1 - \frac1{\ell-1}\,,
\]
which yields the tightness of $(\la_n)$. Indeed,
\[
\cE \setminus \mathcal K_\ell
= \bigcup_{k\ge\ell} \mathcal X_k\,,
\]
where $\mathcal X_k = \bigl\{F\in\cE\colon \max_{S_k}|F|>\mu_k \bigr\}$, and
$ \la_n(\mathcal X_k) = A_{S_n}\bigl\{w\colon \max_{\tau_{-w} S_k} |G|>\mu_k \bigr\}$.
For $1\le n \le n_k-1$ and $w\in S_n$, we have
\[
\max_{\tau_w S_k} |G|  \le \max_{S_n+S_k} |G|  \le \max_{2S_{n_k-1}} |G| < \mu_k
\]
(recall that $k\le n_k-1$), whence, for these $n$s,
$ \bigl\{w\colon \max_{\tau_w S_k} |G|>\mu_k \bigr\} \bigcap S_n = \emptyset$.
On the other hand, for $n\ge n_k$, we have
\[
\la_n(\mathcal X_k) = A_{S_n} \bigl\{w\colon \max_{\tau_w S_k} |G|>\mu_k \bigr\}
\le A_{S_n} \bigl\{w\colon \max_{\tau_w S_k} |G|>M_{n_k} \bigr\} < \frac1{k^2}\,.
\]
Thus,
\[
\la_n (\cE\setminus \mathcal K_l)  \le \sum_{k\ge \ell} \la_n (\mathcal X_\ell)
< \sum_{k\ge \ell} \frac1{k^2} < \frac1{\ell-1}\,,
\]
proving the tightness of $(\la_n)$.

\subsubsection{Translation-invariance of the limiting measure}
Now, let $\la$ be any limiting probability measure for the sequence
$(\la_n)$. Since each measure $\la_n$ is supported by the orbit
$\mathcal O_G$, clearly, $\la$ is supported by the closure of the
orbit $\bar{\mathcal O}_G$.

The measure $\la$ is translation-invariant. This follows from the
fact that for any $n\ge 1$, any $\zeta\in\bC$, and any Borel set $\mathcal X\subset\cE$,
\[
\bigl| \la_n(\tau_\zeta \mathcal X) - \la_n(\mathcal X) \bigr|
\le \frac{A(S_n \triangle \tau_\zeta S_n)}{A(S_n)}
\le \frac{O(|\zeta|)}{L(S_n)}\,,
\]
where $\triangle$ denotes the symmetric difference of sets, and
$L(S_n)$ is the side-length of $S_n$.

\subsubsection{A modification of the limiting measure does not charge the constant functions}
At last, we can specify the measure $\la$ such that it will not charge
the set $\{\rm const\}$ of constant functions. Indeed, following
our assumption~\eqref{eq-KB2} and passing if necessary to some subsequence,
we may assume that a positive limit exists
\[
\lim_{n\to\infty} A_{S_n}\bigl\{w\colon \osc_{\tau_w S}|G|>c \bigr\}
= \alpha>0\,.
\]
This yields that $\la (\cE\setminus \{{\rm const}\}) \ge \alpha >0$.
To see this, let  $U=\{F\in\cE\colon \operatorname{osc}_S |F|<\frac12\,c\}$.
Then $U$ is an open set and $U\supset  \{\text{const}\}$. Hence, it is enough to show
that, for each $n$, $\la_n(U)\le 1-\alpha$. This holds since
\begin{multline*}
\la_n (U) = A_{S_n}\bigl( \bigl\{w\colon \tau_w G \in U  \bigr\} \bigl)
= A_{S_n}\bigl( \bigl\{w\colon \operatorname{osc}_S \tau_w G <\frac{c}2  \bigr\} \bigr)
\\
= A_{S_n}\bigl( \bigl\{w\colon \operatorname{osc}_{\tau_w S} G <\frac{c}2  \bigr\} \bigr)
\le 1-\alpha\,.
\end{multline*}
Then, if needed, we replace $\la$ by its restriction on
$ \cE\setminus \{{\rm const}\} $ and normalize it to make $\la$
the probability measure. This completes the proof of part (i).
\hfill $\Box$

\subsection{Proof of part (ii) of Lemma~\ref{lemma-KB}}

Now, we suppose that condition~\eqref{eq-KB1'} holds, and
assume that the probability measures $(\la_n)$ and $\la$ are the same
as in the proof of part (i). Consider the open set
\[
\mathcal X_k = \bigl\{F\in\cE\colon \max_{S_k} |F| > M_k  \bigr\}\,.
\]
We have
\[
\la_n (\mathcal X_k) =
A_{S_n}\bigl( \bigl\{w\colon \max_{\tau_w S_k} |G| > M_k \bigr\} \bigr)\,,
\]
whence, by~\eqref{eq-KB1'},
\[
\sum_{k\ge 1} \Bigl( \limsup_{n\to\infty} \la_n(\mathcal X_k) \Bigr) < \infty\,.
\]
Furthermore, since the sets $\mathcal X_k$ are open,
\[
\la (\mathcal X_k) \le \limsup_{n\to\infty} \la_n(\mathcal X_k)\,,
\]
so
\[
\sum_{k\ge 1} \la (\mathcal X_k) < \infty\,.
\]
Hence, applying the Borel-Cantelli lemma, we conclude that
\[
\la \Bigl( \bigcap_{\ell\ge 1} \bigcup_{k\ge\ell} \mathcal X_k \Bigr) = 0\,,
\]
which means that $\la$-a.e. $F\in\cE$ does not belong to any $\mathcal X_k$ with
$k\ge k_0(F)$, i.e.,
\[
\limsup_{k\to\infty} \bigl( \max_{S_k} |F| - M_k  \bigr) \le 0\,.
\]
This proves part (ii) and finishes off the proof of Lemma~\ref{lemma-KB}.
\hfill $\Box$

\section{Proof of Theorems~1B,~2B, and~3B}
After the work we have done in Lemmas~\ref{lemma-EF} and~\ref{lemma-KB}, the proofs of these
theorems is rather straightforward.

\subsection{Proof of Theorem~1B}
We take the sequence
\[
\e_j = \frac1{(j+10) \log^3(j+10)}, \quad j\ge 1\,,
\]
put $a_0=1$ and $a_n =  3(1+\e_j)a_{n-1}$ for $n\ge 1$, and
(with some conflict of notation used in Lemma~\ref{lemma-SH})
$S_n'=[-a_n, a_n]^2$.
By $G$ we denote  the corresponding entire function with properties
as in Lemma~\ref{lemma-EF}. We fix a sufficiently large
value of the parameter $B$ as in Lemma~\ref{lemma-EF}
and then will drop dependence on $B$ from our notation.
We claim that
\begin{itemize}
\item
{\em conditions~\eqref{eq-KB1'} and~\eqref{eq-KB2}
of Lemma~\ref{lemma-KB} {\rm (}part {\rm (ii)}{\rm )} are hold for the sequences
$S_n=S_{2^n}'$ and $M_n=\exp\mathcal M(2^{n+1})+1$}.
\end{itemize}

\subsubsection{}
First, we verify convergence of the series
\[
\sum_{k\ge 1} \Bigl(  \limsup_{n\to\infty} A_{S_n} \bigl\{w\colon \max_{\tau_w S_k} |G| > M_k \bigr\} \Bigr)<\infty\,.
\]
For this, we need to bound the relative area
\[
A_{S_{2^n}'} \bigl( \bigl\{w\colon \max_{\tau_w S_{2^k}'} |G|
> \exp\mathcal M(2^{k+1})+1 \bigr\} \bigr)\,.
\]
We note that for $\zeta\in  [-a_{2^{k+1}}+a_{2^k}, a_{2^{k+1}}-a_{2^k}]^2$, the translations $\tau_\zeta S_{2^k}'$
belong to $S_{2^{k+1}}'$. Thus, for $w=w_j(2^{k+1}, 2^n)+\zeta$, $0\le j \le 9^{2^n-2^{k+1}}-1$, we have
\begin{align*}
\max_{\tau_w S_{2^k}'} |G|  &= \max_{\tau_\zeta S_{2^k}'} |\tau_{-w_j(2^{k+1}, 2^n)} G|  \\
&\le \max_{S_{2^{k+1}}'}  |\tau_{-w_j (2^{k+1}, 2^n)} G|  \\
&< \max_{S_{2^{k+1}}'} |G| + \Delta_{2^{k+1}} \\
&< M_k\,.
\end{align*}
The relative area of the set of these $w$s in $S_{2^n}'$ is
\begin{align*}
\frac{9^{2^n-2^{k+1}} (a_{2^{k+1}}-a_{2^k})^2}{a_{2^n}^2} &=
 \frac{9^{2^n-2^{k+1}} (a_{2^{k+1}}-a_{2^k})^2}{9^{2^n-2^{k+1}} a_{2^{k+1}}^2 \prod_{j=2^{k+1}}^{2^n} (1+\e_j)^2} \qquad ({\rm since\ } a_n=3a_{n-1}(1+\e_n))\\ \\
 &= \Bigl( 1-\frac{a_{2^k}}{a_{2^{k+1}}} \Bigr)^2 \cdot
 \Bigl( 1 - (2+o(1)) \sum_{j=2^{k+1}}^{2^n} \e_j \Bigr) \\ \\
 &\ge 1 - 2\, \frac{a_{2^k}}{a_{2^{k+1}}} -   (2+o(1)) \sum_{j\ge 2^{k+1}} \e_j \,, \qquad k\to \infty\,.
\end{align*}
Hence,
\begin{equation}\label{eq:limsup}
\limsup_{n\to\infty} A_{S_n} \bigl\{w\colon \max_{\tau_w S_k} |G| > M_k \bigr\}
\le 2\, \frac{a_{2^k}}{a_{2^{k+1}}}  +  (2+o(1)) \sum_{j\ge 2^{k+1}} \e_j\,.
\end{equation}
At last, for $\ell\to\infty$, $a_\ell = (a+o(1))3^\ell$ with some $a>0$, and
\[
\sum_{j\ge\ell} \e_j = \frac{1+o(1)}{2\log^2\ell}
\qquad (\,{\rm since\ } \e_j=\frac1{(j+10)\log^3(j+10)}\,)\,.
\]
Therefore, the RHS of~\eqref{eq:limsup} is
\[
\lesssim 3^{-2^k} + \frac1{k^2}\,,
\]
which is what we need for condition~\eqref{eq-KB1'}.

\subsubsection{}
To verify condition~\eqref{eq-KB2}, we take $S=[-a_1, a_1]^2$, $\delta<a_1-1$,
and note that for $w=w_j(1, 2^n)+\zeta$ with $0\le j \le 9^{2^n-1}-1$, $|\zeta|<\delta$,
we have $\tau_w S \supset \tau_{w_j(1, 2^n)} [-1, 1]^2$. Therefore,
\begin{align*}
\osc_{\tau_w S} |G| &\ge \osc_{\tau_{w_j(1, 2^n)} [-1, 1]^2} |G| \\
&\ge \osc_{[-1, 1]^2}|G| -\Delta_1 \\
&\ge \osc_{[-1, 1]^2}|z| - \frac13 - \Delta_1 \\
&\ge c > 0
\end{align*}
since $\frac13+\Delta_1 < \frac13+1 <\sqrt{2}=\osc_{[-1, 1]^2} |z|$.
Thus, the set
$\bigl\{w\colon \osc_{\tau_w S}|G|>c \bigr\}$
contains the $\delta$-neighbourhood of the set
$\bigl\{w_j (1, 2^n) \colon 0\le j \le 9^{2^n-1}-1\bigr\} $.
Hence, the relative area of this set in $S_n=S_{2^n}'$ is bounded from below by
\[
\frac{9^{2^n} \pi\delta^2}{a_{2^n}^2} \gtrsim \delta^2 >0\,,
\]
which yields condition~\eqref{eq-KB2}.

\subsubsection{}
At last, applying Lemma~\ref{lemma-KB}, we see that for $\la$-a.e. $F\in\cE$,
\[
\limsup_{[-a_{2^n}, a_{2^n} ]^2} \bigl( |F| - \exp \mathcal M(2^{n+1}) \bigr) \le 1\,.
\]
In our case
$\mathcal M(m) \lesssim \exp(Cm^2 \log^3 m)$, whence
$\mathcal M(2^{n+1})  \le \exp(C  2^{2n} n^3 )$. Then, given $R\ge 10$, we choose $n$ such that
$a_{2^{n-1}} < R \le a_{2^n}$ and  get
\[
\log M_F(R) = \max_{R\, \bD} \log  |F| \le \max_{[-a_{2^n}, a_{2^n} ]^2} \log |F|
\le \exp(C  2^{2n} n^3 )\,.
\]
Furthermore, recalling that $a_m = (a+o(1))3^m$, we see that $2^n \lesssim \log a_{2^{n-1}}$,
whence $2^{2n} n^3 \lesssim (\log R)^2 (\log\log R)^3$, proving Theorem~1B. \hfill $\Box$

\subsection{Proof of Theorem~2B}

Here, we take $\e_j = 3^{-j}$, and let $G$ be the entire function constructed
by using Lemma~\ref{lemma-EF}. Note that
in this case
\[
\mathcal M(n) \le \exp \bigl( C 3^n \bigr) \le \exp \bigl( C a_n \bigr)\,.
\]
Given $R\ge 10$ we choose $n$ so that $a_{n-1}<R \le a_n$ and get
\[
\log M_G (R) \le e^{CR}\,.
\]
Furthermore, given a square $S_k$, for any $n>k$ and any $j\in\bigl\{0, 1, \ldots , 9^{n-k}-1 \bigr\}$, we have
\[
\max_{S_k} \bigl| \tau_{w_j (k, n)} G - G \bigr| < \Delta_k\,.
\]
Given $\e>0$ and a compact set $K\subset\bC$, we choose $k$ so large that $K\subset S_k$ and
\[
\max_{S_k} \bigl| \tau_{w_j (k, n)} G - G\bigr| < \e
\qquad {\rm for\ any\ } n>k\,.
\]
Observing that each square $S\subset \bC$ with the side length $C3^{k+1}$ contains at least one point of the
set
\[
\bigl\{w_j(k, n)\colon 0\le j \le 9^{n-k}-1, n\ge k+1 \bigr\}\,,
\]
we complete the proof of Theorem~2B. \hfill $\Box$

\subsection{Proof of Theorem~3B}

As in the proof of Theorem~2B, we take $\e_j = 3^{-j}$. We denote by $(S_n)$ the corresponding ternary system of squares
and let $G$ be the entire function as in Lemma~\ref{lemma-EF}. We fix $B$ so large that
\[
\max_{S_k} |G|  + \sum_{j\ge 1} \Delta_j < e^{\mathcal M_B(k)}\,,
\]
and drop the parameter $B$ in our notation.  As in the proof of Theorem~1B, a straightforward verification, which we skip,
shows that
conditions~\eqref{eq-KB1} and~\eqref{eq-KB2} of Lemma~\ref{lemma-KB} are satisfied.

As before, we put
\[
\la_n = \frac1{A(S_n)}\, \int_{S_n} \delta_{\tau_w G}\, {\rm d}A(w)\,,
\]
denote by $\la$ any limiting measure and by $(n_i)$ the sequence of indices such
that $\la_{n_i}\to \la$ weakly.

We fix $t$ sufficiently large and choose $k$ so that
$e^{\mathcal M(k-1)} < t \le e^{\mathcal M(k)} $. Then for all $t$'s (except
maybe a countable set of values which we may neglect),
\begin{align*}
\la\bigl\{F\in\cE\colon |F(0)|>t  \bigr\} &=
\lim_{i\to\infty} \la_{n_i} \bigl\{ F\in\cE\colon |F(0)|>t \bigr\} \\ \\
&= \lim_{i\to\infty} A_{S_{n_i}} \bigl\{w\colon |G(w)|>t  \bigr\}\,.
\end{align*}
Since $t>e^{\mathcal M(k-1)}$, we have $|G|<t$ on $S_{k-1}$, as well as on all translations $\tau_{w_j(k-1, n)} S_{k-1}$.
Thus,
\begin{align*}
A_{S_n} \bigl\{w\colon |G(w)|>t  \bigr\} &\le
A_{S_n} \Bigl(  \bigl( \bigcup_{j=0}^{9^{n-k}-1} \tau_{w_j(k-1, n)} S_{k-1} \bigr)^c \Bigr) \\ \\
&= 1 - \frac{9^{n-k} (2a_{k-1})^2}{(2a_n)^2} \\ \\
&= 1 - \prod_{j=k}^n (1+2\e_j)^{-2} \\ \\
&=(2+o(1))\, \sum_{j=k}^n \e_j \\ \\
&\lesssim 3^{-k}\,.
\end{align*}
On the other hand, we have
\[
\log t \le \mathcal M(k) \le e^{C 3^{k}}\,,
\]
whence
\[
3^{-k} \lesssim \frac1{\log\log t}\,,
\]
completing the proof of Theorem~3B. \hfill $\Box$

\bigskip

\medskip

{\noindent L.B.:
 School of Mathematics, Tel Aviv University, Tel Aviv, 69978 Israel,
\newline{\tt levbuh@post.tau.ac.il}

\smallskip
\noindent A.G.:
 School of Mathematics, Tel Aviv University, Tel Aviv, 69978 Israel,
\newline{\tt adiglucksam@gmail.com}

\smallskip
\noindent A.L.:
School of Mathematics, Tel Aviv University, Tel Aviv, 69978 Israel, \newline
\& Chebyshev Laboratory, St. Petersburg State University, 14th Line V.O., 29B,
\newline Saint Petersburg 199178 Russia, {\tt 	log239@yandex.ru}

\smallskip
\noindent M.S.:
 School of Mathematics, Tel Aviv University, Tel Aviv, 69978 Israel,
\newline{\tt sodin@post.tau.ac.il}

}

\end{document}